\documentclass[11pt]{amsart}

\usepackage{amscd,amsthm,amsfonts,amssymb,amsmath,comment,supertabular,mathtools}
\usepackage[only,mapsfrom,Mapsto,Mapsfrom]{stmaryrd}
\usepackage{mathrsfs}
\usepackage[colorlinks=true]{hyperref}
\usepackage{fullpage}
\usepackage{yfonts,dsfont,bbm}
\usepackage[all]{xy}
\usepackage{xcolor}
\usepackage{multirow} 
\usepackage{array}    
\usepackage{graphicx} 

    \newcommand{\BA}{{\mathbb {A}}} 
     
     \newcommand{\BF}{{\mathbb {F}}}

    \newcommand{\BK}{{\mathbb {K}}} 
     
    \newcommand{\BO}{{\mathbb {O}}} 
    \newcommand{\BQ}{{\mathbb {Q}}} \newcommand{\BR}{{\mathbb {R}}}

     \newcommand{\BZ}{{\mathbb {Z}}}
    \newcommand{\bp}{{\mathbbm {p}}}

     \newcommand{\CB}{{\mathcal {B}}}
    \newcommand{\CC}{{\mathcal {C}}} \renewcommand{\CD}{{\mathcal {D}}}

    \newcommand{\CI}{{\mathcal {I}}} 
     \newcommand{\CL}{{\mathcal {L}}}
    \newcommand{\CO}{{\mathcal {O}}} \newcommand{\CP}{{\mathcal {P}}}
     
    \newcommand{\CS}{{\mathcal {S}}}

     \newcommand{\RH}{{\mathrm {H}}}

     \newcommand{\fp}{{\mathfrak{p}}}

     \newcommand{\fP}{{\mathfrak{P}}}

    \newcommand{\Aut}{{\mathrm{Aut}}}

    \newcommand{\End}{{\mathrm{End}}} 
    
    \newcommand{\Fr}{{\mathrm{Fr}}}

    \newcommand{\Gal}{{\mathrm{Gal}}} 
    
    \newcommand{\Hom}{{\mathrm{Hom}}}
    
    \renewcommand{\Im}{{\mathrm{Im}}}

    \newcommand{\Jac}{{\mathrm{Jac}}}\newcommand{\Ker}{{\mathrm{Ker}}}

    \newcommand{\loc}{{\mathrm{loc}}}

    \newcommand{\ord}{{\mathrm{ord}}} \newcommand{\rank}{{\mathrm{rank}}}
     \newcommand{\Pic}{\mathrm{Pic}}

    \renewcommand{\mod}{\ \mathrm{mod}\ }
    
    \newcommand{\ram}{{\mathrm{ram}}}

    \newcommand{\Sel}{{\mathrm{Sel}}}

    \newcommand{\Sym}{{\mathrm{Sym}}}

    \newcommand{\tr}{{\mathrm{tr}}}\newcommand{\tor}{{\mathrm{tor}}}
    
    \newcommand{\ur}{{\mathrm{ur}}}

        \newcommand{\Tr}{\mathrm{Tr}}

    \newcommand{\sD}{\mathscr{D}}

\DeclareFontFamily{U}{wncy}{}
\DeclareFontShape{U}{wncy}{m}{n}{<->wncyr10}{}
\DeclareSymbolFont{mcy}{U}{wncy}{m}{n}
\DeclareMathSymbol{\Sha}{\mathord}{mcy}{"58}

    \newcommand{\ov}{\overline}

    \newcommand{\lra}{\longrightarrow}
    \newcommand{\ra}{\rightarrow}

    \newcommand{\N}{\mathrm{N}}

                    \newcommand{\Res}{\mathrm{Res}}

    \theoremstyle{plain}
    \newtheorem{thm}{Theorem}[section] \newtheorem{coro}[thm]{Corollary}
    \newtheorem{lem}[thm]{Lemma}  \newtheorem{prop}[thm]{Proposition}
     \newtheorem{defn}[thm]{Definition}

\theoremstyle{remark} \newtheorem{remark}{Remark}[section]
\theoremstyle{remark} 
\theoremstyle{remark} 

    \numberwithin{equation}{section}

\newcommand{\rk}{\mathrm{rk}}

\newcommand{\CMS}{{\mathrm {CMS}}} 
\newcommand{\CMU}{{\mathrm {CMU}}} 
\newcommand{\sym}{{\mathrm {Sym}}} 
\newcommand{\uni}{{\mathrm {Uni}}}

\begin{document}
\title[Selmer ranks]{Selmer ranks in twists of CM abelian varieties}
\author[J. Shu]{Jie Shu}
\address{School of Mathematical Sciences, Tongji University, Shanghai 200092,  P. R. China}
\email{shujie@tongji.edu.cn}

\begin{abstract}
We prove the Selmer ranks in certain families of  $p$-th twists of CM abelian varieties obey the symplectic or unitary distributions. As an application, for a prime $p\geq 3$, we obtain that the twisted Fermat curves $X^p+Y^p=\delta$ over a number field containing a primitive $p$-th root of unity are ``largely" unsolvable as $\delta$ varies. We also discuss the rank growth in cyclic extensions of prime degree for CM abelian varieties. 

\end{abstract}


\thanks{}

\maketitle


\section{Introduction}\label{Intro-Selmer}
\subsection{Twisted Fermat curves}
 Let $p\geq 3$ be a prime number. Fermat's last theorem, proved by A. Wiles \cite{Wiles95}, states that the Fermat equation $X^p+Y^p=1$ has no $\BQ$-solutions with $XY\neq 0$.  A. Diaconu and Y. Tian \cite{DT05} prove, over certain totally real field $F$, there exist infinitely many $p$-th power-free $\delta\in F^\times$ such that the twisted Fermat equation $X^p+Y^p=\delta$ has no $F$-solutions.  There is also a bulk of works concerning the twisted Fermat cubics over $\BQ$, e.g. \cite{Lieman94, CST17, ABS22, KS24}.

Let $F$ be an algebraic number field containing the $p$-th roots of unity.  Under certain ``fan-structure" stratification on $F^\times/{F^\times}^p$, we prove 
\begin{thm}\label{Main}
The proportion of $\delta\in F^\times/{F^\times}^p$ such that the twisted Fermat equation $X^p+Y^p=\delta$ has no solutions in $F$ is at least $\sD_p^\Sym(0):=\prod_{i=1}^\infty(1+p^{-i})^{-1}$. 
\end{thm}

\begin{proof}
 Let  $C_{\delta}$ be the smooth curve over $F$ defined by the affine equation $y^p=x( \delta-x)$, and let $J_{\delta}$ be its Jacobian variety. 
 Denote
\[\Pi_{1}=\{\delta\in F^\times/F^{\times p }: \rk_{\BZ}\  J_{\delta}(F)=0.\}\]
By the same argument using Northcott's finiteness theorem as in \cite{DT05}, one can prove for all but finitely many $\delta\in \Pi_1$, $X^p+Y^p=\delta$ has no $F$-rational solutions.  The theorem follows from Corollary \ref{Fermat-rk0} which provides a lower bound for the set $\Pi_1$ along the ``fan-structure".
\end{proof}

By Proposition \ref{distribution}, $\sD_p^\Sym(0)>1-1/(p-1)$. If $p$ varies and we vary $F$ accordingly with the $p$-th roots of unity contained in $F$, then we immediately have the following corollary.
\begin{coro}
As $p$ goes to $\infty$, the probability of the equation $X^p+Y^p=\delta$ has an $F$-solution with $\delta\in F^\times/F^{\times p}$ tends to zero. 
\end{coro}

 The proportion or probability under the ``fan-structure" stratification will be made precise in the next subsection. The proof of Theorem \ref{Main} is based on our general results on the distribution of Selmer groups in twists of CM abelian varieties which we describe in the following.

\subsection{Selmer ranks in twists of CM abelian varieties}
Let $F$ be an algebraic number field with $G_F$ its absolute Galois group. Let $\BK$ be a CM field with complex conjugation $\dag$ and denote by $\BK_0$ its maximal totally real subfield. Let $\CO_\BK$ be the ring of integers of $\BK$ and $\BO$ an order of $\BK$ which is stable under $\dag$. Let $A/F$ be an abelian variety and denote  $\End_F^0(A)=\End_F(A)\otimes \BQ$. Assume 
\begin{equation}\tag{{\bf CM}}\label{CM}
[\BK:\BQ]=2\dim A \text{ and there is an embedding $\BK\hookrightarrow \End^0_F(A)$ with $\BO=\BK\cap \End_F(A)$.}
\end{equation}
Let $A^\vee$ be the dual abelian variety of $A$. The CM action of $A$ induces an embedding $\BO\hookrightarrow \End_{F}(A^\vee)$ via $a\mapsto a^\vee$. Fix a symmetric isogeny $\lambda:A\ra A^\vee$ that is $\dag$-sesquilinear, i.e. $\lambda a^\dag=a^\vee\lambda$ for all $a\in \BO$. The $\dag$-sesquilinearity of $\lambda$ exactly means that the Rosati involution induced by $\lambda$ restricts to the involution $\dag$ of $\BO$.

Let $p$ be a rational prime. Let $\bp$ be a prime ideal of $\BO$ above $p$ such that $\bp^\dag=\bp$.  Let $\omega$ be a primitive $p$-th root of unity and $\mu_p=\langle \omega\rangle$. Denote $\CO=\BO[\omega]$, $K=\BK(\omega)$ and $\fp=(\bp,1-\omega)$ the unique prime ideal of $\CO$ above $\bp$. Assume
\begin{itemize}
\item[({\bf H0})] If $\omega\notin \BO$, then $\CO$ is a finitely generated projective $\BO$-module. For example one may take $\BO=\CO_\BK$.
\item[({\bf H1})] $p\nmid [\CO_{K}:\CO]\deg(\lambda)$. Then $\fp$ resp. $\bp$ is an invertible prime ideal of $\CO$ resp. $\BO$ and denote $k=\BO/\bp=\CO/\fp=\BF_q$ the residue field.
\item[({\bf H2})]  $T=A[\bp]$ is a simple $\BF_p[G_F]$-module.
\item[({\bf H3})]   One of the following holds:
\begin{itemize}
\item[$(\CMS)$] $\bp$ is ramified in $\BK/\BK_0$ and $k=\BF_p$.
\item[$(\CMU)$] $\bp$ is inert in $\BK/\BK_0$ and $k=\BF_{p^2}$.
\end{itemize}
\end{itemize}
\begin{remark}
\begin{itemize}
\item[(1)] Under $({\bf H0})$, if $p\nmid [\CO_K:\CO]$, then one can prove  $p\nmid [\CO_\BK:\BO]$.
\item[(2)] If $(\CMS)$ holds, $A[\bp]\cong \BO/\bp=\BF_p$, and thus $({\bf H_2})$ holds automatically. 
\end{itemize}
\end{remark}

Let $\CC(F)=\Hom(G_F,\mu_{p})$ be the group of continuous characters of order dividing $p$. For any $\chi\in \CC(F)$, denote $\CO_\chi$ to be the $G_F$-module whose underlying group is $\CO$ and $g\cdot x=\chi(g)x$ for $g\in G_F$ and $x\in \CO$. The Serre tensor $ A_\chi=\CO_\chi\otimes_\BO A$ is an abelian variety over $F$ with CM by $\CO$. Through the Kummer $\fp$-descent map $\kappa_{A_\chi,\fp}$ (cf. (\ref{kummer-map})), the $\fp$-Selmer group
\[\Sel_\fp(A_\chi)=\Ker\left(\RH^1(F,A_\chi[\fp])\ra\prod_{v}\RH^1(F_v,A_\chi)\right)\]
fits into the Kummer exact sequence
\[0\ra A_\chi(F)\otimes_\CO k\xrightarrow{\kappa_{A_\chi,\fp}} \Sel_\fp(A_\chi)\lra \Sha(A_\chi)[\fp]\ra 0.\]

As explained in \cite{Shu-quadratic}, via the duality theorems,  the CM-action and the symmetric isogeny $\lambda$ induce a quadratic structure on the restricted product 
\[\RH^1(\BA_F,A[\bp])={\prod_v}'\left(\RH^1(F_v, A[\bp]),\RH^1_\ur(F_v, A[\bp])\right),\]
 which becomes a metabolic symplectic or unitary $k$-space according to $(\CMS)$ or $(\CMU)$ respectively. By Theorem \ref{twist-av}, $A[\bp]\cong A_\chi[\fp]$. This identification induces an identification of the $k$-spaces $\RH^1(\BA_F,A[\bp])\cong \RH^1(\BA_F,A_\chi[\fp])$. As $\chi$ varies in $\CC(F)$, the Selmer group $\Sel_\fp(A_\chi)$ runs as the intersections of the global cohomology and various Kummer images, which are all maximal isotropic subspaces of $\RH^1(\BA_F,A[\bp])$. Therefore, one may expect to describe the behavior of these Selmer groups in terms of the symplectic or unitary distribution, i.e. the rank distribution of the intersections of a fixed maximal isotropic subspace with a random one in a metabolic infinite-dimensional symplectic or unitary $k$-space (cf. \cite[\S 1.1]{Shu-quadratic}). If $(\CMS)$ holds, take $\epsilon=\sym=1$; Otherwise take $\epsilon=\uni=1/2$. 
 
 \begin{defn}
 The symplectic $(\epsilon=1)$ or unitary $(\epsilon=1/2)$ distribution $\sD_q^\epsilon$ is given as,  for non-negative integers $r$,
\[\sD_q^\epsilon(r)=\prod_{i\geq 0}(1+q^{-i-\epsilon})^{-1} \prod_{i=1}^r\frac{q^{1-\epsilon}}{q^i-1}.\]
\end{defn}

Let $\CB$ be a finite set of places of $F$ containing all archimedean places, places above $p$ and places where $A$ has degenerate reductions. Let $\CO_{\CB}$ be the ring of $\CB$-integers of $F$. Suppose $\CB$ is large enough so that $\Pic(\CO_{\CB})=0$. Denote $d=\ord_\fp(1-\omega)$ and, for $0\leq i\leq d$, put $F_i=F(A[\bp^i])$. For every $k\geq -1$ and $X>0$, let $\Delta_k(X)$ be the ``fan-structure" of collections of characters of $G_F$ of order dividing $p$ in the spirit of Klagsbrun-Mazur-Rubin \cite{KMR14}. The ``fan-structue" depends on $\CB$ and $A[\bp^d]$, and we have $\CC(F)=\bigcup_{k,X}\Delta_k(X)$ (cf. Definition \ref{fan}, \ref{charset}, \ref{stratification0} and \ref{stratification}). Denote $\Delta(k,X)=\bigcup_{i=-1}^k\Delta_i(X)$.
\begin{thm}\label{main}
Assume either $F=F_{d-1}\neq F_d$ with $d\geq 2$ or $F_1=F_d$. Then
\[\lim_{k\ra \infty}\lim_{X\ra \infty} \frac{\left|\{\chi\in \Delta(k,X):\dim_k \left(\Sel_\fp(A_\chi)/\kappa_{A_\chi,\fp}(A_\chi(F)_\tor)\right)=r\}\right|}{\left|\Delta(k,X)\right|}=\sD_q^\epsilon(r),\]
where $\epsilon=1$ resp. $1/2$ if $(\CMS)$ resp. $(\CMU)$ holds.
\end{thm}
Theorem \ref{main} is a special case of Theorem \ref{main-body} and Corollary \ref{main2-body} and their variants Theorem \ref{main-body1} and Corollary \ref{main2-body1}. The ``fan-structrue" stratification $\CC(F)=\bigcup_{k,X} \Delta_k(X)$ is in fact a disjoint union for $k$ and $\bigcup_{X>0} \Delta_k(X)$ consists of characters ramified exactly at $k$ (resp. $k+1$) ``critical" places if $F\neq F_d$ (resp. $F=F_d$). Thus $\bigcup_{X>0} \Delta(k,X)$ consists of characters for which the number of ``critical" ramified places is  $\leq k$ (resp. $\leq k+1$) and 
\begin{equation}\label{ordering}
\lim_{k\ra \infty}\lim_{X\ra \infty} \Delta(k,X)=\CC(F).
\end{equation}
In the following, the densities, averages or proportions all mean the limits along (\ref{ordering}). 
\begin{coro}\label{main2}
Assume either $F=F_{d-1}\neq F_d$ with $d\geq 2$ or $F_1=F_d$.
\begin{itemize} 
\item[(1)] The average size of the $\fp$-Selmer groups modulo torsions is $1+q^{1-\epsilon}$. More precisely,
\[\lim_{k\ra\infty}\lim_{X\ra\infty}\frac{\sum_{\chi\in \Delta(k,X)}\left|\Sel_{\fp}(A_\chi)/\kappa_{A_\chi,\fp}(A_\chi(F)_\tor)\right|}{\left|\Delta(k,X)\right|}=1+q^{1-\epsilon}.\]
\item[(2)] The twists of rank zero have density $\geq \sD_q^\epsilon(0)>1-\frac{q^{1-\epsilon}}{q-1}.$
\item[(3)] The average Mordell-Weil rank $\leq [K:\BQ]\sum_{i=0}^\infty \frac{1}{1+q^{i+\epsilon}}.$
\end{itemize}
\end{coro}
\begin{proof}
(1) and (2) follows from Theorem \ref{main} and Proposition \ref{distribution} immediately. Note the localization $\CO_{(\fp)}$ is a discrete valuation ring. For the torsion-free $\CO$-module $M_\chi=A_\chi(F)/A_\chi(F)_\tor$, $M_{\chi}\otimes_\CO\CO_{(p)}$ is a free $\CO_{(\fp)}$-module of finite rank. Therefore $\dim_k M_\chi\otimes_\CO k=\dim_K M_\chi\otimes_\CO K$, and thus 
\[\frac{1}{[K:\BQ]}\rank_\BZ A_\chi(F)=\dim_k M_\chi \otimes_\CO k\leq \dim_k \Sel_{\fp}(A_\chi)/\kappa_{A_\chi,\fp}(A_\chi(F)_\tor).\]
On the other hand, by Theorem \ref{main} and Proposition \ref{distribution}, $\Sel_{\fp}(A_\chi)/\kappa_{A_\chi,\fp}(A_\chi(F)_\tor)$  has average $k$-dimension $\sum_{i=0}^\infty \frac{1}{1+q^{i+\epsilon}}$. Thus (3) follows as desired.  
\end{proof}

\subsection{Rank growth in cyclic extensions}
Let $A/F$ be an abelian variety. For a fixed prime $p$, let $\Xi$ denote the set of cyclic extensions over $F$ of degree $p$. As the question raised by Mazur, Rubin and Larsen (\cite[\S2]{Mazur83}, \cite[\S C]{MRL18}), it is natural to ask how ``large" is the set 
\[\CS_p(A/F)=\{L/F\in \Xi: \rank_\BZ A(L)>\rank_\BZ A(F)\}.\]
Conditional on the Birch and Swinnerton-Dyer conjecture, $\CS_p(A/F)$ is closely related to the collection of characters $\chi$ of $F$ of order $p$ such that the $L$-function $L(s,A,\chi)$ of the abelian variety twisted by $\chi$ vanishes at the central point $s=1$. There is a good deal of statistical heuristics of such zeros, particularly in the case where $E$ is an elliptic curve over $\BQ$. These heuristics predict the magnitude of $\CS_p(E/\BQ)$ decreases as $p$ grows and is finite for $p\geq 7$ (see \cite{DFK07,MR23} for precise statements). In contrast, over general number fields there can be a large supply of cyclic extensions $L/F$ in which the rank grows. For example, the rank frequently grows in extensions of dihedral type, as can be detected for parity reasons, or nontrivial Heegner points (cf. \cite{Cornut02,Vatsal03,MR07}).

Let $\pi:\CC(F)\backslash \{1\}\ra \Xi$ be the map sending $\chi$ to the fixed field of $\Ker(\chi)$. Then $\pi$ is surjective and each fibre consists of a Galois orbit of characters of order $p$. Since the characters in a single fibre have the same ramification places,  the ``fan-structure" ordering transfers well to $\Xi$ through $\pi$ and we write $\Xi=\bigcup_{k,X}\Xi_k(X)$. If we denote $A_1=\CO_1\otimes_{\BO}A$ for the trivial character $\chi=1$, for $L\in \Xi$, $\rank_\BZ A_1(L)$ grows if and only if $\rank_\BZ A(L)$ grows. Moreover,  $A_\chi$ is the twist of $A_1$ by the character $\chi^{-1}$ (cf. Definition \ref{st} and the following paragraph). For an abelian group $M$, denote $M_\BQ=M\otimes_\BZ \BQ$. Thus for each $L\in \Xi$, we have the decomposition
\begin{equation}\label{eigen-decomposition}
A_1(L)_\BQ=A_1(F)_\BQ\bigoplus \left(\bigoplus_{\chi\in \pi^{-1}(L)}A_\chi(F)_\BQ\right).
\end{equation}
In particular, we can estimate the rank growth of $A_1$, or equivalently $A$, by the distribution of the Selmer groups of $A_\chi$.
\begin{thm} \label{growth}
Let the notations and conditions be as in Theorem \ref{main}. There is a positive proportion of cyclic extensions $L\in \Xi$  such that $ \rank_\BZ A(L)=\rank_\BZ A(F)$. Moreover, if $p=2$, the positive proportion is at least $\sD_q^{\epsilon}(0)$.
\end{thm}

If $p=2$, the map $\pi$ is a bijection and the density of trivial (i.e. minimal) Selmer groups transfers to give a lower bound. For $p\geq 3$, in proving Theorem \ref{growth}, the lower bound obtained for the positive proportion is small.  In fact, we use the worst estimation pretending that for each $L\in \Xi$ there are at least $p-2$ characters of $\pi^{-1}(L)$ for which the twists have trivial Selmer groups. In general we expect, in spirit of the heuristics of \cite{DFK07,MR23}, the proportion of cyclic extensions in $\Xi$ without rank growth for $A$ is large, at least for large primes $p$. This is true for some subfamilies.  Let $\Xi^1\subset \Xi$ be the subset of cyclic extensions which split completely over places of $\CB$. Then $\CC^1(F)=\pi^{-1}(\Xi^1)$ consists exactly of characters restricting to trivial local characters at the places of $\CB$. 
\begin{thm}\label{growth1}
Let the notations and conditions be as in Theorem \ref{main} and suppose $\omega\notin \BO_\bp$.  The proportion of cyclic extensions $L\in \Xi^1$  such that $ \rank_\BZ A(L)=\rank_\BZ A(F)$ is at least $\sD_q^\epsilon(0)$.
\end{thm}
For the above theorem, the key is that for $L\in \Xi^1$, the Selmer groups associated to characters of $\pi^{-1}(L)$ are all isomorphic. Generally, we don't expect Selmer groups associated to Galois conjugacy characters to be isomorphic (cf. Proposition \ref{local-condition2}.(3)). For example for twists of Jacobian varieties of Fermat curves, this is generally not  true as can be detected by the root-number formulae or explicit dimension formulae for Selmer groups in \cite{Shu21jnt, Shu-Selmer}.

\begin{thm}\label{growth2}
Let the notations and conditions be as in Theorem \ref{main}. Suppose $(\CMS)$ holds with $\BO=\CO_{\BK}$. Suppose 
\begin{itemize}
\item if $p=2$, then the discriminant $d_{\BK_0}$ is odd;
\item if $p\neq 2$, then $p$ is prime to the relative discriminant $d_{\BK/M}$  and $[\BQ(\omega): M]$ is odd where $M=\BQ(\omega)\cap \BK$. 
\end{itemize}
Assume $\Sha(A_\chi)[\fp^\infty]$ is finite for all $\chi\in \CC(F)$.
\begin{itemize}
\item[(1)] The proportion of  $\chi\in \CC(F)$ with $\rank_\CO A_\chi(F)$ odd is $\frac{1-\beta^\sym}{2}$ where $\beta^\sym=\prod_{i=1}^\infty\frac{1-p^{-i}}{1+p^{-i}}$.
\item[(2)] The proportion of cyclic extensions $L\in \Xi$  such that $ \rank_\BZ A(L)>\rank_\BZ A(F)$ is at least $\frac{1-\beta^\sym}{2}$.
\end{itemize}
\end{thm}

Theorem \ref{growth2} applies to quadratic twists of elliptic curves with CM by the ring of integers of an imaginary quadratic field in which the prime $2$ ramifies. Theorem \ref{growth2} also applies to $p$-th twists of elliptic curves with CM by $\BZ\left[\frac{-1+\sqrt{-p}}{2}\right]$ with $p\equiv 3\mod 4$ and the Jacobian varieties of Fermat curves (cf. Section \ref{examples}). The proofs of Theorem \ref{growth}-\ref{growth2} are delayed into the finial section.

\subsection{Some general remarks}
Generally, the character group $\CC(F)$ is divided into a finite union of subfamilies $\CC^\gamma(F)$ where $\gamma\in \Omega_1=\prod_{v\in \CB}\CC(F_v)$. The subfamily $\CC^\gamma(F)$ consists of characters restricting to $\gamma$ over $\CB$. Accordingly, the twists $\{{A_\chi:\chi\in \CC(F)}\}$ are divided into several subfamilies. The Selmer groups $\Sel_\fp(A_\chi)$ with $\chi\in \CC^\gamma(F)$ always contain two kinds of constant contributors: One arises from the torsions $\kappa_{A_\chi,\fp}(A_\chi(F)_\tor)$; Another is a subspace $\Sha_{\CP_1}(\gamma)$ that only depends on $\gamma$ (cf. Definition \ref{constant}). The distribution of $\Sel_\fp(A_\chi)$ with $\chi\in \CC^\gamma(F)$ is described in term of $\sD_q^\epsilon$, shifted by the dimension of $(\kappa_{A_\chi,\fp}(A_\chi(F)_\tor)+\Sha_{\CP_1}(\gamma))$ (cf. Lemma \ref{torsion-gamma}, Theorem \ref{main-body}). Particularly, under the hypothesis of Theorem \ref{main},  for any $\gamma\in \Omega_1$, we have $\Sha_{\CP_1}(\gamma)\subset \kappa_{A_\chi,\fp}(A_\chi(F)_\tor)$ and thus the shift only depends on torsions. One may refer to Lemma \ref{torsion-gamma} for a discussion of these two subspaces. 

The two cases where $\omega\in \BO$ or not are treated very differently at some critical points. If $\omega\in \BO$, then $A_\chi$ are all twists of $A$. Then much geometric information of $A$ can be carried onto $A_\chi$ through the twisting map. For example, the quadratic structures on $\RH^1(F_v, A[\bp])$ and $\RH^1(F_v,A_\chi[\fp])$ coincide under the twisting isomorphism. If $\omega\notin \BO$, then $A_\chi$ are twists of $A_1=\CO_1\otimes_\BO A$. But here some geometric obstruction occurs. Generally there is no symmetric isogeny (or polarization) on $A_\chi$ of degree prime to $p$ whose Rostai involution will restrict to the complex conjugation of $\CO$ (cf. \cite{Howe01}). In such situations, we can't equip $\RH^1(F_v,A_\chi[\fp])$ with (non-degenerate) quadratic structures and this results in different approaches to prove that local Kummer images are maximal isotropic subspaces (cf. Proposition \ref{id-WP} and \ref{isotropic}). There are also different patterns in the equidistribution of local Kummer images (cf. Proposition \ref{local-condition3} and \ref{local-condition2}). Moreover, in Appendix \ref{stc}, we consider certain situations where we can construct symmetric isogenies on $A_1$, and hence on all $A_\chi$, of degree prime to $p$ with Rosati involution restricting to the complex conjugation of $\CO$ and these constructions will be used in proving Theorem \ref{growth2}. 

Indeed, we prove all the results in the setting of $p^n$-th twists. We note, for higher twists, the distributions remain the same, up to a shift depending on $\gamma$, because the local Kummer images twisted by a character $\chi$ of order $p^n$ are essentially determined by the character $\chi^{p^{n-1}}$ of order $p$ (cf. Proposition  \ref{local-condition3}-\ref{equi}). The Markovian machinery used to derive the distributions is due to Klagsbrun, Mazur and Rubin \cite{KMR14} where they developed the machinery to treat orthogonal cohomological spaces. Instead, for the situation of CM abelian varieties, we are dealing with symplectic and unitary cohomological spaces. 

We end this introduction by including the following numerical table for the readers' convenience which should be read in companion with Proposition \ref{distribution}.
\begin{table}[h]
\begin{tabular}{|ll|l|l|l|l|l|l|}
\hline
\multicolumn{2}{|r|}{$p\quad\quad\quad\quad$}                     & 2 & 3&5&7& 11&13\\ \hline
\multicolumn{1}{|l|}{\multirow{2}{*}{$\quad\sD_q^\epsilon(0)$}} & $\sym$ & 0.4194& 0.6390&0.7933&0.8545&0.9084&0.9226 \\ \cline{2-8} 
\multicolumn{1}{|l|}{}                  &  $\uni$& 0.5686 & 0.7198&0.8264&0.8724&0.9159&0.9281 \\ \hline
\multicolumn{1}{|l|}{\multirow{2}{*}{$\sum_{\text{odd}}\sD_q^\epsilon(r)$}} & $\sym$  &0.4394 &0.3210 &0.1984 &0.1424&0.0908&0.0768 \\ \cline{2-8} 
\multicolumn{1}{|l|}{}                  & $\uni$ &0.3807 &0.2699 &0.1721 &0.1272&0.0839&0.0718 \\ \hline
\multicolumn{1}{|l|}{\multirow{2}{*}{$\sum r\sD_q^\epsilon(r)$}} & $\sym$  &  0.7644&0.4040&0.2150 &0.1483&0.0923& 0.0778\\ \cline{2-8} 
\multicolumn{1}{|l|}{}                  & $\uni$ & 0.4850 &0.2903  &0.1749&0.1279&0.0840&0.0718\\ \hline
\end{tabular}
\caption{Probabilities of zero, odd variables and expectations for the symplectic and unitary distributions.}
\end{table}
\subsection*{Acknowledgements}
We thank Ye Tian, Lian Duan and Zikang Wang for helpful communications.

\section{Examples and applications}\label{examples}
\subsection{Elliptic curves}
For a CM elliptic curve $E$, we always take $\lambda$ to be the canonical theta polarization. 

\subsubsection{Congruent number elliptic curves.}
Let $i=\sqrt{-1}$, $F=K=\BQ(i)$, $\CO_K=\BZ[i]$ and $\fp=(1-i)$. Fix $b\in F^\times/F^{\times 4}$. Let $E_b:y^2=x^3-b x$ be an elliptic curve over $F$ and $E_b$ has CM by $\CO_K$: $i(x,y)=(-x,i  y)$. Then  $E_b[\fp]=\{O, (0,0)\}$ and $E_b[2]=\{\CO,(0,0),(\sqrt{b},0),(-\sqrt{b},0)\}$. 

\noindent{$\bullet$ \em Quadratic twists.} Let $\omega=-1$ and then $d=2$ and $F=F_1$. The element $b\in F^{\times 2}$ if and only if $F_1=F_2$. Under Hilbert's Satz 90, we identify 
\[\CC(F)\cong F^\times/F^{\times 2}, \quad \chi_\delta \mapsfrom \delta, \]
where $\chi_\delta(\sigma)=(\sqrt{\delta})^{\sigma-1}$. Then $E_{b\delta^2}$ is the quadratic twist of $E_b$ by $\chi_\delta$.  By Theorem \ref{main}, for $r\geq 1$,
\[\lim_{k\ra \infty}\lim_{X\ra \infty} \frac{\left|\{\delta \in \Delta(k,X):\dim_{\BF_2}\Sel_{(1-i)}(E_{b\delta^2})=r\}\right|}{\left|\Delta(k,X)\right|}=\sD_2^\Sym(r-1).\]

\noindent{$\bullet$ \em Quartic twists.} Let  $\omega=i$ and then $d=1$ and $F=F_1$. Under Hilbert's Satz 90, we identify 
\[\Delta(F)\cong \{\delta\in F^\times/F^{\times 4}: \text{either $2\nmid \ord_v(\delta)$ or $4\mid \ord_v(\delta)$ for all $v\nmid 2b$}\}, \quad \chi_\delta \mapsfrom \delta, \]
where $\chi_\delta(\sigma)=(\sqrt[4]{\delta})^{\sigma-1}$ (cf. Definition \ref{charset}). Then $E_{b\delta}$ is the quartic twist of $E_b$ by $\chi_\delta$.  By Lemma \ref{torsion-gamma}, for any $\gamma\in \Omega_1$, $\Sha_{\CP_1}(\gamma)=0$ and Corollary \ref{main2-body} implies, for $r\geq 1$,
\[\lim_{k\ra \infty}\lim_{X\ra \infty} \frac{\left|\{\delta \in \Delta(k,X):\dim_{\BF_2}\Sel_{(1-i)}(E_{b\delta})=r\}\right|}{\left|\Delta(k,X)\right|}=
\sD_2^\sym(r-1).\]

\subsubsection{Fermat cubic curves}
Let $\omega=\frac{-1+\sqrt{-3}}{2}$ be a primitive $3$-rd root of unity and $F=K=\BQ(\omega)$. Let $\CO_K=\BZ[\omega]$ and $\fp=(\sqrt{-3})$. Fix $b\in F^\times/F^{\times 6}$ and  let $E_b:y^2=x^3+b$ be an elliptic curve over $F$. Then $E_b$ has CM by $\CO_K=\BZ[\omega]$: $\omega(x,y)=(\omega x,y)$.   The torsion group $E_b[\sqrt{-3}]=\{O, (0,\pm\sqrt{b})\}$, $d=1$ and $F_1=F(\sqrt{b})$. Thus $F_1=F$ if and only if $b\in F^{\times 2}$.
Under Hilbert's Satz 90,
\[ \CC(F)\cong F^\times/F^{\times 3}, \quad  \chi_\delta\mapsfrom \delta,\]
where $\chi_\delta(\sigma)=(\sqrt[3]{\delta})^{\sigma-1}$ for all $\sigma\in G_F$. Then $E_{b\delta^2}:y^2=x^3+b\delta^2$ is the cubic twist of $E_b$ by $\chi_\delta$. Then 
\[\lim_{k\ra \infty}\lim_{X\ra \infty} \frac{\left|\{\delta \in \Delta(k,X):\dim_{\BF_3}\Sel_{(\sqrt{-3})}(E_{b\delta^2})=r\}\right|}{\left|\Delta(k,X)\right|}=\left\{\begin{aligned}
&\sD_3^\Sym(r),\quad&&\text{if  $b\notin F^{\times 2}$};\\
&\sD_3^\sym(r-1),\quad&&\text{if $b\in F^{\times 2}$}.
\end{aligned}\right.\]

\subsubsection{$p$-th twists of CM elliptic curves}\label{p-twists}
Let $\BK$ be an imaginary quadratic field and $F$ its Hilbert class field. Let $E$ be an elliptic curve with CM by $\BO=\CO_\BK$ over $F$. 
Let $p$ be a prime and let $K=\BK(\omega)$ where $\omega$ is a primitive  $p$-th root of unity. 

\noindent {$\bullet$ ($\CMS$) \em cases.} Assume $\BK=\BQ(\sqrt{-p})$ where $p\equiv 3\mod 4$. Then $\BK\subset K=\BQ(\omega)$ and $\CO_K=\CO_\BK[\omega]$. Let $\bp=(\sqrt{-p})\subset \CO_\BK$ and $\fp=(1-\omega)\subset \CO_K$. The abelian variety $A_\chi=\CO_{K,\chi}\otimes_{\CO_\BK} E$ has CM by $\CO_K$ over $F$ and $A_\chi[\fp]\cong E[\bp]$. Since $E$ has CM by $\CO_\BK$ over $F$, $E[\bp]$ has $\BF_p$-dimension one and thus is a simple $\BF_p[G_F]$-module. We are in the $(\CMS)$ case with $\CO_K/\fp\cong \CO_\BK/\bp\cong \BF_{p}$ and $d=1$. For $\CC(F)=\Hom(G_F,\mu_p)=\bigcup_{k,X}\Delta_k(X)$, Theorem \ref{main} implies, 
\[\lim_{k\ra \infty}\lim_{X\ra \infty} \frac{\left|\{\chi\in \Delta(k,X):\dim_{\BF_{p}} \Sel_{(1-\omega)}(A_\chi)=r\}\right|}{\left|\Delta(k,X)\right|}=\left\{\begin{aligned}
&\sD_p^\Sym(r),\quad&&\text{if  $E[\bp](F)=0$};\\
&\sD_p^\sym(r-1),\quad&&\text{if $E[\bp](F)\neq 0$}.
\end{aligned}\right.\]

\noindent {$\bullet$ ($\CMU$) \em cases.} Let $p$ be a prime inert in $\BK$ and $\bp=(p)\subset \CO_\BK$. Then $\CO_K=\CO_\BK[\omega]$ and the prime $\fp=(1-\omega)\subset \CO_K$ above $p$ is inert in $K/K_0$.  The abelian variety $A_\chi=\CO_{K,\chi}\otimes_{\CO_\BK} E$ has CM by $\CO_K$ over $F$ and $A_\chi[\fp]\cong E[p]$. By \cite[Theorem 6.18]{BC20}, $E[p]$ is a simple $\BF_p[G_F]$-module. We are in the $(\CMU)$ case with $\CO_K/\fp\cong \CO_\BK/\bp\cong \BF_{p^2}$ and $d=1$ (cf. Appendix \ref{STC}).  Theorem \ref{main} implies, for any $r\geq 0$, 

\[\lim_{k\ra \infty}\lim_{X\ra \infty} \frac{\left|\{\chi\in \Delta(k,X):\dim_{\BF_{p^2}}\Sel_{(1-\omega)}(A_\chi)=r\}\right|}{\left|\Delta(k,X)\right|}=
\sD_{p^2}^\uni(r).\]

\subsection{Fermat curves}\label{Fermat}
Let $p\geq 3$ be a prime number. Let $F$ be an algebraic number field containing $\mu_{p}=\langle \omega\rangle$ and let $K=\BQ(\omega)$ be the $p$-th cyclotomic field.  For $\delta\in F^\times/F^{\times p}$, the twisted $p$-th  Fermat curve $W_{\delta}$ over $F$ is defined by the equation
\[X^p+Y^p=\delta.\] 
For any integer $1\leq t\leq p-2$, let  $C_{\delta,t}$ be the smooth curve over $F$ of genus $g=(p-1)/2$ defined by the affine equation $y^p=x^t( \delta-x).$ There is a quotient map given as
\[\varphi_{\delta,t}: W_\delta\lra C_{\delta,t},\quad (X,Y)\mapsto (X^p,X^tY).\]
The Jacobian variety $J_{\delta,t}=\Jac(C_{\delta,t})$ has dimension $g=(p-1)/2$.  The action of  $\omega$ on $C_{\delta,t}$ by $\omega(x,y)=(x,\omega y)$ induces the complex multiplication $\iota: \CO_K\ra \End_F(J_{\delta,t})$. 
\begin{lem}
Let $\lambda:J_{1,t}\ra J_{1,t}^\vee$ be the canonical principal polarization induced by the theta divisor. Then the Rosati involution associated to $\lambda$ restricts to the complex conjugation of $\iota(\CO_K)$. 
\end{lem}
\begin{proof}
Simply denote $C=C_{1,t}$ and $J=J_{1,t}$.  We first unravel the definition of $\lambda$. Recall the theta divisor $\Theta$ is the image of the natural map
\[C^{(g-1)}\lra \Pic^{g-1}_{C/F},\quad (P_1,\cdots,P_{g-1})\mapsto \left[\sum_{i=0}^{g-1} P_i\right].\]
Let $\infty=[0:1:0]$ be the point at infinity on $C$ and $D=[(g-1)\infty]\in  \Pic^{g-1}_{C/F}(F)$. We have the translation $t_D: \Pic_{C/F}^0\xrightarrow{\sim}\Pic^{g-1}_{C/F}$ and denote $\Theta_D=t^*_D \Theta$. Then the polarization $\lambda$ can be explicitly given as
\[\lambda(x)= (\Theta_D+x)-\Theta_D\in J^\vee.\]
In order to check the stability of $\CO_K$, it suffices to show that $\lambda \omega^{-1}=\omega^\vee \lambda$. Since $\omega$ acts trivially on $\Theta$ and $\infty$, $\omega$ fixes $\Theta_D$. Then for any $x\in J$,
\[\omega^\vee\lambda(x)=\omega^{-1}(\Theta_D+x)-\Theta_D)=(\Theta_D+\omega^{-1}x)-\Theta_D)=\lambda \omega^{-1}(x).\]
\end{proof}

Through Hilbert's Satz 90, we identify $F^\times/F^{\times p}\cong\CC(F)$ via  $\delta \mapsto \chi_\delta$ where$\chi_\delta(\sigma)=(\sqrt[p]{\delta})^{\sigma-1}$.
Then  $J_{\delta,t}$ is the twist of $J_{1,t}$ by the character $\chi_{\delta^{t+1}}^{-1}$. Let $\fp=(1-\omega)$ be the unique prime ideal of $\CO_K$ above $p$.  The torsion $J_{1,t}[\fp]=\langle [(0,0)-\infty] \rangle$.  Thus $d=1$ and $F=F_1$. Applying Theorem \ref{main}, we have 
\begin{thm} \label{Fermat-density}
Fix $1\leq t\leq p-2$. For $r\geq 1$,
\[\lim_{k\ra \infty}\lim_{X\ra \infty} \frac{\left|\{\delta\in \Delta(k,X):\dim_{\BF_p}(\Sel_\fp(J_{\delta,t}))=r\}\right|}{\left|\Delta(k,X)\right|}=
\sD_p^\sym(r-1).\]
\end{thm}

Denote
\[\Pi_{t}=\{\delta\in F^\times/F^{\times p }: \rk_{\BZ} \ J_{\delta,t}(F)=0.\}\]

\begin{coro}\label{Fermat-rk0}
For all $1\leq t\leq p-2$,  the set $\Pi_{t}$ has density at least $\sD_p^\Sym(0)$.
\end{coro}
\begin{proof}
The corollary follows from Theorem \ref{Fermat-density} and the exact sequence
\[0\ra J_{\delta,t}(F)\otimes_\CO k\ra \Sel_\fp(J_{\delta,t})\ra \Sha(J_{\delta,t})[\fp]\ra 0, \]
together with the fact $\dim_{k} J_{\delta,t}(F)_\tor\otimes k=1$.
\end{proof}

\section{Twists of abelian varieties}
Since we will deal with the general $p^n$-th twists for a fixed integer $n\geq 1$, we assume $\omega$ is a primitive $p^n$-th root of unity and $\mu_{p^n}=\langle\omega\rangle$. To avoid some technical complexity, if $n\geq 2$, we assume $\mu_{p^n}\subset F$.  If $A$ is absolutely simple with CM type $(\BK,\Phi)$ and $\BK/\BQ$ is abelian, then the reflex type is $(\BK,\Phi^{-1})$, and by \cite[Chap. II, Proposition 30]{Shimura-CM}, $F$ always contains $\BK$. Thus in this case, if $\BK$ contains $\mu_{p^n}$, then $F$ contains $\mu_{p^n}$ automatically.

Assume the data $(A,\lambda,\BK,\BO,p,\bp)$  satisfies (\ref{CM}) and $({\bf H0})$-$({\bf H3})$. 
Let  $\CC(F)=\Hom(G_F,\mu_{p^n})$ be the group of continuous characters of order dividing $p^n$.  For $\chi\in \CC(F)$, let $\CO_\chi$ be the $G_F$-module whose underlying group is $\CO$ and $g\cdot x=\chi(g)x$ for $g\in G_F$ and $x\in \CO$.
\begin{defn}\label{st}
For any $\chi\in \CC(F)$, let $A_\chi=\CO_\chi\otimes_\BO A$ be the Serre tensor of $A$ by $\CO_\chi$ $($cf. \cite[Definition 1.1]{MRS07}, \cite[\S III-1.3]{Serre94}$)$. 
\end{defn}
In view of $({\bf H0})$, $\CO$ is a finitely generated projective $\BO$-module. For $\chi=1$,  $A_1=\CO_1\otimes_\BO A$ is an abelian variety over $F$ (cf. \cite[Theorem 7.2]{Conrad04}) and  $A_\chi$ is the twist of $A_1$ by the character $\chi^{-1}$ in the usual sense. Recall $\fp$ is an invertible prime ideal of $\CO$ above $p$ and $d=\ord_\fp(1-\omega)$.

\begin{thm}\label{twist-av}
Suppose $\chi\in \CC(F)$.
\begin{itemize}
\item[(1)] There is an embedding $\CO\hookrightarrow \End_F(A_\chi)$ and $A_\chi$ has CM by $K$ over $F$.
\item[(2)] The Tate module $$T_{\fp}(A_\chi)\cong\CO_{\fp,\chi}\otimes_{\BO_\bp} T_\bp(A)$$ is a free $\CO_{\fp}$-module of rank one.
\item[(3)] For $1\leq i\leq d$, the twisting induces $G_F$-isomorphisms $A_\chi[\fp^i]\cong A[\bp^i]$.
\end{itemize}
\end{thm}
\begin{proof}
We have $\CO=\End_{\CO[G_F]}(\CO_\chi)$ and the embedding $\BO\hookrightarrow \End_F(A)$. By \cite[Proposition 1.6]{MRS07}, there is an embedding 
\[\CO=\End_{\CO[G_F]}(\CO_\chi)\otimes_{\BO}\BO\hookrightarrow \End_F(A_\chi).\]
Since 
\[[K:\BQ]=[K:\BK][\BK:\BQ]=2[K:\BK]\dim A=2\dim A_\chi,\]
$A_\chi$ has CM by $K$ over $F$. This proves (1).

Since $p$ is prime to the conductor of $\BO$, (2) follows from \cite[Theorem 2.2]{MRS07} and \cite[\S4, Remark 2]{Serre-Tate68}. As for (3), for $1\leq i\leq d$, $(1-\omega)\subset \fp^i$ and thus $ \CO/\fp^i\cong\BO/\bp^i$. Then we have
\begin{eqnarray*}
A_\chi[\fp^i]&\cong&T_\fp(A_\chi)/\fp^i T_\fp(A_\chi)\cong ( \CO_{\fp,\chi}/\fp^i)\otimes_{\BO} T_\bp(A)\\
&\cong & ( \CO/\fp^i)\otimes_{\BO} T_\bp(A)\cong(\BO/\bp^i)\otimes_\BO T_\bp(A)\cong A[\bp^i].
\end{eqnarray*}
\end{proof}

If $\omega\in \BO$, then $A_1=A$ and $A_\chi$ is a twist of $A$. There exists a twisting isomorphism $\phi_{A,\chi}: A\ra A_\chi$ over $\ov{F}$ such that $\chi^{-1}(\sigma)=\phi_{A,\chi}^{-1}\phi_{A,\chi}^{\sigma}\in \Aut_F(A)$. The existence of $\phi_{A,\chi}$ is unique up to an isomorphism over $F$. The assignment $f\mapsto f_\chi:=\phi_{A,\chi}f\phi_{A,\chi}^{-1}$ clearly induces an isomorphism of the algebras $\End_{\ov{F}}^0(A)$ and $\End_{\ov{F}}^0(A_\chi)$. We embed $\BK\hookrightarrow \End_{\ov{F}}(A_\chi)^0$ through this isomorphism. By \cite[Proposition 2.1]{Howe01}, $\BK$ indeed embeds into $\End_F^0(A_\chi)$ and $\BO=\BK\cap \End_F^0(A_\chi)$. In particular, the abelian variety $A_\chi$ satisfies $(\ref{CM})$.
The map $\lambda_\chi:={(\phi_{A,\chi}^{-1})}^\vee\lambda \phi_{A,\chi}^{-1}$ descends to a  $\dag$-sesquilinear symmetric isogney on $A_\chi$. Moreover, if $\lambda$ is a polarization, so is $\lambda_\chi$ (cf.  \cite[Proposition 2.2]{Howe01}). 

For any positive integer $m$, let $e_{A,m}:A[m]\times A^\vee[m]\ra \mu_m$ be the Weil pairing of $A$. Composing the symmetric isogeny $\lambda$, we get an alternating pairing $e_{A,m}^\lambda:A[m]\times A[m]\ra \mu_m$ so that $e_{A,m}^\lambda(x,y)=e_{A,m}(x,\lambda y)$. 

\begin{prop}\label{id-WP}
Suppose $\omega\in \BO$. The twisting map $\phi_{A,\chi}$ identifies the pairings  $e_{A,m}^\lambda$ and $e_{A_\chi,m}^{\lambda_\chi}$.
\end{prop}
\begin{proof}
Simply denote $\phi=\phi_{A,\chi}$. By \cite[Proposition 11.21]{Edixhoven-AV}, 
\[e_{A_\chi,m}^{\lambda_\chi}(\phi x,\phi y)=e_{A_\chi,m}(\phi x,(\phi^{-1})^\vee\lambda y)=e_{A,m}(x, \lambda y)=e_{A,m}^\lambda(x,y)\]
as desired.
\end{proof}

\section{Twisted incoherent Selmer groups}
Let $\Sigma$ be the set of all  places of $F$. For each $v\in \Sigma$, let $F_v$ denote the completion of $F$ at $v$ and  let $\CO_{v}$ be its valuation ring if $v$ is finite. Let $\ov{F}$ resp. $\ov{F_v}$ denote a fixed algebraic closure of $F$ resp. $F_v$ with absolute Galois group $G_F=\Gal(\ov{F}/F)$ resp. $G_v=\Gal(\ov{F_v}/F_v)$.
\subsection{Local Tate pairings and quadratic spaces}
The Weil pairing 
\[E_p:T_p(A)\times T_p(A^\vee)\ra \BZ_p(1)\]
is a $G_F$-equivariant perfect bilinear pairing satisfying $E(a x,y)=E(x,a^\vee y)$ for all $a\in \BO$. Composing $\lambda$, we obtain a non-degenerate alternating pairing
\begin{equation}\label{Weil1}
E_p^\lambda:T_p(A)\times T_p(A)\ra \BZ_p(1),\quad E_p^\lambda(x,y)=E(x,\lambda(y)), 
\end{equation}
which is $\dag$-adjoint, i.e. $E_p^\lambda(ax,y)=E_p^\lambda(x,a^\dag y)$ for all $a\in \BO$.
The pairing $E_p^\lambda$ factors over the decomposition
\[T_pA=\bigoplus_{\fP\mid p} T_{\fP} A\]
where $\fP$ runs over the prime ideals of $\BO$ above $p$. Noting $\bp^\dag=\bp$, the $\bp$-component 
\[E^\lambda_\bp: T_\bp(A)\times T_{\bp}(A)\ra \BZ_p(1)\]
induces an isomorphism
\[T_{\bp}(A)\cong \Hom_{\BZ_{p}}(T_\bp(A),\BZ_p(1)).\]

Let $\BO_0=\BO\cap \BK_0$ and $\bp_0=\bp\cap \BO_0$. Since $p\nmid [\CO_\BK:\BO]$, both $\BO_\bp$ and $\BO_{0,\bp_0}$ are rings of integers of $\BK_\bp$ and $\BK_{0,\bp_0}$ respectively. Fix a generator $d\in \partial^{-1}_{\BO_\bp/\BZ_p}$ of the inverse different so that $d^\dag/d=-1$ resp. $+1$ in the $(\CMS)$ resp. $(\CMU)$ case. As in \cite[Lemma A.3 and Proposition A.4]{MR07}, composing the trace map 
\[\BO_\bp\ra \BZ_p,\quad a\mapsto \Tr_{\BK_\bp/\BQ_p}(da) \]
induces the following isomorphisms
\[\Hom_{\BO_{\bp}}(T_\bp(A),\BO_\bp(1))\cong \Hom_{\BZ_p}(T_\bp(A),\BZ_p(1)) \cong T_{\bp}(A),\]
which in turn defines, by adjunction, a $\dag$-sesquilinear $G_F$-equivariant perfect pairing
\begin{equation}\label{Weil2}
\Theta^\lambda_\bp: T_\bp(A)\times T_{\bp} A\ra \BO_\bp(1),
\end{equation}
characterized by
\begin{equation}\label{characterization}
E^\lambda_\bp(ax,y)=E^\lambda_\bp(x,a^\dag y)=\Tr_{\BK_\bp/\BQ_p}(da \Theta_\bp^\lambda( x,y) ), \quad \forall \ x,y\in T_\bp(A), a\in \BO_\bp.
\end{equation}

Denote $T=A[\bp]$ and, reducing modulo $\bp$, $\Theta_\bp^\lambda$ yields a $G_F$-equivariant   perfect  pairing 
\[\theta_\bp^\lambda: T\times T\ra k(1).\]
For each $v\in \Sigma$, the cup product with respect to the pairing $\theta_\bp^\lambda$ induces the local Tate pairing
\begin{equation}\label{tate-pairing}
h_v=\cup_{\theta_\bp^\lambda}:\RH^1(F_v,T)\times\RH^1(F_v,T)\lra \RH^2(F_v,k(1))\xhookrightarrow[]{\mathrm{Inv}_v}  k.
\end{equation}
The last inclusion is the Hasse invariant map  which is an  isomorphism for nonarchimedean places.
As a special case of \cite[Proposition 4.4 and Theorem 5.4]{Shu-quadratic}, we have
\begin{thm}\label{Hermitian}
\begin{itemize}
\item[(1)] In the $(\CMS)$ case, $(\RH^1(F_v,T), h_v)$ is a metabolic symplectic  $k$-space.
\item[(2)] In the $(\CMU)$ case, $(\RH^1(F_v,T), h_v)$ is a metabolic unitary $k$-space.
\end{itemize}
\end{thm}

For a subspace $X\subset \RH^1(F_v,T)$, denote by $X^\perp$ the orthogonal complement of $X$ under $h_v$. Then $X$ is called maximal isotropic if $X= X^\perp$. Denote
\[\CI_v=\{\text{maximal isotropic subspaces of $(\RH^1(F_v,T),h_v)$}\}.\]
If $T$ is unramified at $v$ and  $\RH^1_\ur(F_v,T)=\Ker(\RH^1(F_v,T)\ra \RH^1(F_v^\ur, T))$ is the subgroup of unramified cohomologies, define
\[\CI^\ram_v=\{X\in \CI_v: X\cap \RH^1_\ur(F_v,T)=0\}.\]

\begin{prop}\label{ur-cohom}
Suppose $v\in \Sigma$ is a finite place not dividing $p$ and $T$ is unramified at $v$. 
\begin{itemize}
\item[(1)] $\RH^1_\ur(F_v,T)$ is a maximal isotropic subspace in $\RH^1(F_v,T)$.
\item[(2)] $\dim_{k} \RH^1(F_v,T)=2\dim_{k} T^{G_v}$.
\item[(3)] If $\dim_kT^{G_v}=1$, then $(\RH^1(F_v,T),h_v)$ is a hyperbolic plane and $|\CI^\ram_v|=p$.
\end{itemize}
\end{prop}
\begin{proof}
By \cite[Theorem 2.4]{Tateduality}, $\RH^1_\ur(F_v,T)$ is maximal isotropic with respect to $h_v$ and
\[\dim_{k} T^{G_v}=\dim_{k}\RH^1_\ur(F_v,T).\]
Consequently,
\[\dim_{k} \RH^1(F_v,T)=2\dim_{k}T^{G_v}.\]

If $T^{G_v}=T$, then, by the proceeding assertion, $\RH^1(F_v,T)$ has dimension $2$. By Theorem \ref{Hermitian}, $\RH^1(F_v,T)$ is metabolic, and hence a hyperbolic plane. By \cite[Lemma 2.17]{Shu-quadratic}, there are $1+p$ (maximal) isotropic lines. By Assertion (1), $\RH^1_\ur(F_v,T)$ is one of them and thus $|\CI^\ram_v|=p$ as desired.
\end{proof}

\subsection{Equidistribution for local Kummer conditions}
For $\chi\in  \CC(F_v)$, we describe the $\fp$-descent on $A_\chi$ as follows (see \cite[Appendix A]{GP12, Shu-quadratic}). For $F$-schemes $T$, the functor $T\mapsto \fp^{-1}\otimes_\CO A_\chi(T)=\Hom_\CO(\fp,A(T))$ is represented by an abelian variety over $F$, denoted as $\fp^{-1}\otimes_\CO A_\chi$. The inclusion $\CO\hookrightarrow \fp^{-1}$ induces an isogeny $A_\chi\ra \fp^{-1}\otimes_{\CO}A_\chi$ with kernel $A_\chi[\fp]$, i.e. there is the exact sequence 
\[0\ra A_\chi[\fp]\ra A_\chi\xrightarrow{} \fp^{-1}\otimes_{\CO} A_\chi\ra 0.\]
For each $v\in \Sigma$, taking Galois cohomology, the boundary map gives rise to the Kummer $\fp$-descent map
\begin{equation}\label{kummer-map}
\kappa_{A_\chi, \fp,v}:\fp^{-1}/\CO\otimes_{\CO}A_\chi(F_v)\hookrightarrow \RH^1(F_v,A_\chi[\fp]).
\end{equation}
Through the identification $A_\chi[\fp]\cong A[\bp]$ in Theorem \ref{twist-av}.(3), define 
\[\CL_{A,v,\chi}:=\Im\left( A_\chi(F_v)/\fp A_\chi(F_v)\xrightarrow{\kappa_{A_\chi,\fp,v}}\RH^1(F_v, A_\chi[\fp])\xrightarrow{\cong} \RH^1(F_v, T)\right).\]
\begin{prop}\label{isotropic}
For any $\chi\in \CC(F_v)$, $\CL_{A,v,\chi}$ is a maximal isotropic subpace in $\RH^1(F_v,T)$.
\end{prop}
\begin{proof}
We first give a proof in the case $\omega\in \BO$. While identifing $A_\chi[\fp]$ with $A[\bp]$, by Proposition \ref{id-WP} the twisting map $\phi_{A,\chi}$ also identifies the pairings $e_{A_\chi, p^m}^{\lambda_\chi}$ and $e_{A,p^m}^\lambda$. Thus  $\phi_{A,\chi}$ identifies $\RH^1(F_v, A_\chi[\fp])\xrightarrow{\cong} \RH^1(F_v, T)$ together with their quadratic structures. The proposition follows from the fact that
the Kummer image $\Im(\kappa_{A_\chi,\fp,v})$ is maximal isotropic in $\RH^1(F_v, A_\chi[\fp])$ which is a consequence of the classical Tate dualities (cf. \cite[Proposition 5.5 and 5.6]{Shu-quadratic}).

For the case $\omega\notin \BO$, extending the involution $\dag$ to be the complex conjugation on $\CO$, there is generally no $\dag$-sesquilinear symmetric isogenies of $A_\chi$ of degree prime to $p$. So we can't equip quadratic structures on $\RH^1(F_v,A_\chi[\fp])$. In this case, we follow an idea of Mazur and Rubin \cite[Proposition A.7]{MR07} using a generalization of Tate dualities due to Bloch-Kato. Of course this treats the case $\omega\in \BO$ uniformly. Define a $G_F$-equivariant perfect pairing 
\[f:\CO_{\chi,\fp}\times \CO_{\chi,\fp}\ra \CO_\fp,\quad f(x,y)=xy^\dag.\]
 Through the isomorphism 
\[T_\fp(A_\chi)\cong {\CO_\chi}\otimes_{\BO} T_\bp(A),\] 
 the Weil pairing (\ref{Weil2}) and $f$ define a pairing
\[\langle\ ,\ \rangle_{\CO_\fp}:=f\otimes_{\BO_\bp} \Theta_\bp^\lambda: T_{\fp}(A_\chi)\times T_{\fp}(A_\chi)\lra \CO_\fp(1).\]
This pairing is perfect, $G_F$-equivariant and $\dag$-sesquilinear (with respect to $\CO_\fp$ in the second variable). By Pontryagin duality, $\langle\ ,\ \rangle_{\CO_\fp}$ induces a pairing
\[\langle\ ,\ \rangle_\fp:\left(T_{\fp}(A_\chi)\otimes_{\BZ_p}\BQ_p/\BZ_p\right)[\fp]\times T_{\fp}(A_\chi)/\fp T_{\fp}(A_\chi)\lra \fp^{-1}\CO_\fp/\CO_\fp(1).\]
The identifications (which depend on choices of uniformizers of $\CO_\fp$ and $\BO_{\bp}$)
\begin{equation}\label{id}
\left\{\begin{aligned}
&\left(T_{\fp}(A_\chi)\otimes_{\BZ_\bp} \BQ_p/\BZ_p\right)[\fp]=\fp^{-1}/\CO\otimes_\CO T_\fp(A_\chi)\cong \CO/\fp\otimes_\CO T_\fp(A_\chi)\cong A_\chi[\fp],\\
&A_\chi[\fp]\cong \CO/\fp\otimes_\CO T_\fp(A_\chi)\cong \BO/\bp\otimes_\BO T_\bp(A)\cong A[\bp]
\end{aligned}\right.
\end{equation}
transform the pairing $\langle\ ,\ \rangle_{\fp}$ to $\theta_\bp^\lambda$.

By a generalization of Tate duality due to Bloch-Kato \cite[Proposition 3.8, Example 3.11]{BK90}, 
for every place $v$ of $F$, the pairing $\langle\ ,\ \rangle_{\CO_\fp}$ induces a perfect cup-product pairing
\[\cup_{\langle,\rangle_{\CO_\fp}}:\RH^1(F_v,T_{\fp}(A_\chi)\otimes_{\BZ_p} \BQ_p/\BZ_p)\times \RH^1(F_v, T_{\fp}(A_\chi))\lra K_\fp/\CO_\fp,\]
under which the image of $A_\chi(F_v)\ra \RH^1(F_v,T_{\fp}(A_\chi))$ and the image of $A_\chi(F_v)\otimes_{\BZ_p} \BQ_p/\BZ_p\ra \RH^1(F_v,T_{\fp}(A_\chi)\otimes \BQ_p/\BZ_p)$ are orthogonal complements of each other. In fact, let $\langle\ ,\ \rangle_{\BZ_p}$ be the unique perfect $\dag$-adjoint $\BZ_p$-bilinear associated to $\langle\ , \rangle_{\CO_\fp}$ via \cite[Corollary B.3]{Shu-quadratic} such that $\langle x,y\rangle_{\BZ_p}=\Tr_{K_\fp/\BQ_p}(d\langle x,y\rangle_{\CO_\fp}).$ By \cite[Proposition 3.8, Example 3.11]{BK90}, the Kummer images of $A_\chi(F_v)$ and $A_\chi(F_v)\otimes_{\BZ_p} \BQ_p/\BZ_p$ are orthogonal complements of each other under $\cup_{\langle,\rangle_{\BZ_p}}$. Since 
\[\cup_{\langle,\rangle_{\BZ_p}}(\xi,\eta)=\Tr_{K_\fp/\BQ_p}(d\cup_{\langle,\rangle_{\CO_\fp}}(\xi,\eta)),\]
by \cite[Proposition B.4]{Shu-quadratic}, these two Kummer images are also orthogonal complements of each other under $\cup_{\langle,\rangle_{\CO_\fp}}$.

Again by Pontryagin duality,  $\cup_{\langle,\rangle_{\CO_\fp}}$ induces the perfect pairing
\[\cup_{\langle,\rangle_\fp}:\RH^1\left(F_v,\left(T_{\fp}(A_\chi)\otimes_{\BZ_p} \BQ_p/\BZ_p\right)[\fp]\right)\times \RH^1\left(F_v, \CO/\fp\otimes_{\CO}T_{\fp}(A_\chi)\right)\lra \fp^{-1}\CO_\fp/\CO_\fp\cong k.\]
The identifications (\ref{id}) transform $\cup_{\langle,\rangle_{\fp}}$ into the local Tate pairing $h_v=\cup_{\theta_\bp^\lambda}$ in (\ref{tate-pairing}). Thus through Kummer maps (for example \cite[Example 3.11]{BK90}), both the image of $A_\chi(F_v)$ under
\[\RH^1(F_v,T_{\fp}(A_\chi))\ra \RH^1(F_v,T_{\fp}(A_\chi)/\fp T_{\fp}(A_\chi))\xrightarrow{\sim}\RH^1(F_v,T)\]
and the inverse image of $A_\chi(F_v)\otimes \BQ_p/\BZ_p$ under
\[\RH^1(F_v,T)\xrightarrow{\sim}\RH^1(F_v,\left(T_{\fp}(A_\chi)\otimes \BQ_p/\BZ_p\right)[\fp] )\ra \RH^1(F_v,T_{\fp}(A_\chi)\otimes \BQ_p/\BZ_p)\]
are equal to $\CL_{A,v,\chi}$ and are orthogonal complements of each other under $h_v$ as desired.

\end{proof}

Let $\CB\subset \Sigma$ be a finite set of places of $F$ containing all archimedean places, all places above $p$ and all places where $A$ has bad reductions. 
\begin{defn}
For $0\leq i\leq d$, define $F_i=F(A[\bp^i])$ and
\[\CP_0=\{v\notin \CB: \text{$v$ is not split in $F_d$}\},\quad  \CP_1=\{v\notin\CB: \text{$v$ is split in $F_d/F$}\}.\]
\end{defn}
Then $\Sigma=\CB\bigcup \CP_0\bigcup \CP_1$. If $F=F_d$, then $\CP_{0}$ is empty; Otherwise $\CP_0$ has positive density among places of $F$ by Chebotarev's density theorem.

\begin{prop}\label{local-condition1}
\begin{itemize}
\item[(1)] For $v\in \CP_0$ and $\chi\in\CC(F_v)$, 
\[\CL_{A,v,\chi}=\RH^1_\ur(F_v, T)\subset \RH^1(F_v,T).\]
\item[(2)] For $v\notin \CB$ and $\chi\in \CC(F_v)$ an unramified character, 
\[\CL_{A,v,\chi}=\RH^1_\ur(F_v, T)\subset \RH^1(F_v,T).\]
\end{itemize}
\end{prop}
\begin{proof}
Let  $v\in \CP_0$ and $\chi\in \CC(F_v)$. There exists an $0\leq i\leq d-1$ such that $F_v=F_v(A[\bp^i])$ but $F_v\neq F_v(A[\bp^{i+1}])$.  By Theorem \ref{twist-av}, $A[\bp^d]\cong A_\chi[\fp^d]$ as $G_v$-modules. Since $v\nmid p\infty$, $A_\chi(F_v)=A_\chi[\fp^i]\times B$ where $B$ is a profinite group such that $\fp B=B$. Note $A$ has good reduction at $v$.  By the criterion of N\'eron-Ogg-Shafarevich \cite{Serre-Tate68},  the field $E=F_v(A_\chi[\fp^{i+1}])=F_v(A[\bp^{i+1}])$ is unramified over $F_v$.  Then $A_\chi(F_v)\subset \fp A_\chi(E)$, and thus, 
 \[\CL_{A,v,\chi}\subset \RH^1(\Gal(E/F_v), T)\subset \RH^1_\ur(F_v, T).\] 
 The equlity follows by Proposition \ref{ur-cohom} and \ref{isotropic} that both sides have dimension one.
 
 Next suppose $v\notin \CB$ and $\chi\in \CC(F_v)$ is unramified. Let $L_\chi$ be the field over $F_v$ fixed by $\Ker(\chi)$. Then $A$ and hence $A_1$ have good reduction at $v$ and $A_1\cong A_\chi$ over $L_\chi$. Since $L_\chi$ is unramified over $F_v$, it follows from \cite[Corollary 4]{Serre-Tate68} that $A_\chi$ has good reduction over $F_v$.  Since $|\CO/\fp|$ is invertible in $\CO_v$, by \cite[Lemma 6 ]{GP12}, $\CL_{A,v,\chi}=\RH^1_\ur(F_v,T)$.
\end{proof}
\begin{remark}
The places in $\CP_{0}$ are so-called ``free places'' as in \cite{KMR14} for twisting the Selmer groups. Indeed, by the above proposition, the choice of the local characters at places in $\CP_{0}$ doesn't affect the twisted incoherent Selmer groups (see Definition \ref{Selmer group}, \ref{incoherent Selmer}).
\end{remark}

\begin{defn}
 Let $v$ be a finite place of $F$. A character $\chi\in \CC(F_v)$ is totally ramified if  $L_\chi/F_v$ is a totally ramified extension where $L_\chi\subset \ov{F_v}$ is the fixed field by $\Ker(\chi)$.  
 \end{defn}
Define
\[ \begin{aligned}
 &\CC_{n}(F_v)=\{\chi\in \CC(F_v): \text{$\chi$ is of  order $p^n$} \}, \\
 &\CC_{\leq n-1}(F_v)=\CC(F_v)\backslash \CC_n(F_v), \text{ and }\\
 &\CC_{n}^\tr(F_v)=\{\chi\in \CC_n(F_v): \text{$\chi$ is totally ramified} \}.
 \end{aligned}\]
If $v\in \CP_1$, then $\mu_{p^n}\subset F_v$ and we have
 \[|\CC_n(F_v)|=p^{2n-2}(p^2-1),\quad |\CC_{\leq n-1}(F_v)|=p^{2n-2}\text{ and }|\CC_n^\tr(F_v)|=p^{2n-1}(p-1).\]
\begin{lem}\label{character}
 Let $v$ be a finite place of $F$, and $\chi\in \CC_n(F_v)$. We also view $\chi$ as a character of $F_v^\times$ via local class field theory.
 \begin{itemize}
 \item[(1)] $\chi\in \CC_n^\tr(F_v)$ if and only if $\chi(\CO_v^\times)=\mu_{p^n}$, if and only if $\chi^{p^{n-1}}$ is ramified.
 \item[(2)] If $\chi\notin \CC_n^\tr(F_v)$ has order $p^n$, then $\chi=\chi_0\mu$ for some unramified character  $\chi_0\in \CC_n(F_v)$ and a character $\mu\in\CC_{\leq n-1}(F_v)$.
 \end{itemize}
\end{lem}
\begin{proof}
The character $\chi\in \CC_{n}(F_v)$ is totally ramified if and only if $L_{\chi}$ is of degree $p^n$ over $F_v$ and  contains no unramified proper subextensions of $F_v$. Equivalently, by local class field theory, for any $0\neq i\in \BZ/p^{n}\BZ$, the restriction of $\chi^i$ on $\CO_v^\times$ is notrivial which is also equivalent to $\chi(\CO_v^\times)=\mu_{p^n}$. Also for $\chi\in \CC_n(F_v)$, $\chi(\CO_v^\times)=\mu_{p^n}$ if and only if $\chi^{p^{n-1}}(\CO_v^\times)=\mu_p$, i.e. $\chi^{p^{n-1}}$ is ramified. The first assertion follows.

If $\chi\notin \CC_n^\tr(F_v)$, by Assertion (1), $\chi^{p^{n-1}}$ is unramified of order $p$. There exists  an unramified character $\chi_0\in \CC_n(F_v)$ such that $\chi_0^{p^{n-1}}=\chi^{p^{n-1}}$, and thus $\chi/\chi_0\in \CC_{\leq n-1}(F_v)$ as desired.
\end{proof}

By Theorem \ref{twist-av},  the Tate module $T_\bp(A)$ is a free $\BO_\bp$-module of rank one. Thus $G_F$ acts on $T_\bp(A)$ through a character, say,
\begin{equation}\label{CM-character}
\theta: G_v\ra \Aut_\BO(T_\bp(A))=\BO_\bp^\times.
\end{equation}
 For each $v$, denote $\theta_v$ to be the restriction of $\theta$ to $G_v$.

\begin{lem}\label{gen}
Suppose $v\in \CP_1$ and $\chi\in \CC_{n}^\tr(F_v)$.
\begin{itemize}
\item[(1)] If $\omega\notin\BO_\bp$, then $d=1$ and
\[A_\chi[\fp^2]\cong {\CO_\chi}/\fp^2\otimes_{\BO}A[\bp]\text{ and }A_\chi[\fp^{2}](F_v)=A_\chi[\fp].\]
\item[(2)] If $\omega\in\BO_\bp$, then 
\[A_\chi[\fp^{d+1}]\cong{\CO_\chi}\otimes_\BO A[\bp^{d+1}]\text{ and }A_\chi[\fp^{d+1}](F_v)=A_\chi[\fp^{d}].\]
\end{itemize}
\end{lem}
\begin{proof}
(1) First suppose $\omega\notin \BO_\bp$. By Lemma \ref{cyc-ext}, $\fp$ is ramified over $\bp$ and $d=1$. Then 
\[\CO_\chi/\fp^2=\CO_\chi/\fp^2\otimes_\BO \BO=\CO_\chi/\fp^2\otimes_\BO \BO/\bp,\]
and thus, by Theorem \ref{twist-av},
\begin{eqnarray*}
\label{def-field1}
A_\chi[\fp^{2}]&\cong&{\CO_\chi}/{\fp^2}\otimes_\BO T_\bp(A)={\CO_\chi}/\fp^2\otimes_{\BO}\BO/\bp\otimes_\BO T_\bp(A)\\
&=&{\CO_\chi}/\fp^2\otimes_{\BO}A[\bp].\nonumber
\end{eqnarray*}
Since $v\in \CP_1$, $G_v$ acts trivially on $A[\bp]$. Then $\gamma\in G_v$ acts on $A_\chi[\fp]$ by $\gamma\otimes 1$. Let $g\in \Gal(L_\chi/F_v)$ be a generator of the Galois group. Since $\chi$ is totally ramified of oder $p^n$, $\chi(g)$ is a primitive $p^n$-th root of unity by Lemma \ref{character}. Then
\begin{eqnarray*}
\label{def-field}
A_\chi[\fp^{2}](F_v)&=&\left(A_\chi[\fp^2]\right)^{g-1}=(\CO_\chi/\fp^2\otimes_\BO A[\bp])^{g\otimes 1-1}\\
&=&A_\chi[\fp^{2},1-\chi(g)]=A_\chi[\fp].
\end{eqnarray*}

\noindent (2) Suppose $\omega\in \BO_\bp$. Then $\CO_\fp=\BO_\bp$ and $\fp\CO_\fp=\bp\BO_\bp$. By Theorem \ref{twist-av}, we have 
\[A_\chi[\fp^{d+1}]\cong \CO_{\fp,\chi}\otimes_{\BO_\bp} A[\bp^{d+1}]\] and thus
$G_v$ acts on $A_\chi[\fp^{d+1}]$ through the character $\chi\theta_v$ on $A[\bp^{d+1}]$. Since $v\notin \CB$, $A$ has good reduction at $v$ and hence $\theta_v$ is an unramified character.  Note $\chi$ is totally ramified of order $p^n$ and the complex multiplication injects $\mu_{p^n}/\mu_{p^{n-1}}$ into $\Aut_\BO(A[\bp^{d+1}])$ by the definition of $d$.  By considering the effect of an element in the inertia group $I_v\subset G_v$ which restricts to a generator $\Gal(L_\chi/F_v)$, we see  the action of $G_v$ on $A[\bp^{d+1}]$ through the character $\chi\theta_v$ is nontrivial. Therefore $G_v$ acts on $A_\chi[\fp^{d+1}]$ non-trivially and thus $A_\chi[\fp^{d+1}](F_v)=A_\chi[\fp^{d}]$ as desired.

\end{proof}

\begin{prop}\label{local-condition3}
Suppose $\omega \notin \BO_\bp$, $v\in \CP_1$,  and $\chi\in \CC_n(F_v)$. Let $E_\chi\subset L_\chi$ be the unique sub-extension over $F_v$ of degree $p$. Then
\[\CL_{A,v,\chi}=\Hom(\Gal(E_\chi/F_v),T)\subset \RH^1(F_v,T).\]
\end{prop}
\begin{proof}
By Lemma \ref{gen}.(1), $A_\chi[\fp^2](F_v)=A_\chi[\fp]$. Thus, since $v\nmid p\infty$, $A_\chi(F_v)=A_\chi[\fp]\times B$ where $B$ is a profinite group such that $\fp B=B$. Since $v\in \CP_1$, $G_v$ acts trivially on $A[\bp]$ and thus, again by Lemma \ref{gen}.(1), $G_{L_\chi}$ acts trivially on $A_\chi[\fp^2]$, and hence $A_\chi[\fp^2]\subset A_\chi(L_\chi)$. Then $A_\chi(F_v)\subset\fp A_\chi(L_\chi)$ and
\[\CL_{A,v,\chi}\subset \RH^1\left(\Gal(L_{\chi}/F_v),T\right)=\Hom\left(\Gal(L_\chi/F_v),T\right).\] 
By Proposition \ref{isotropic}, both sides of the above inclusion have dimension one, and thus we have equalities
\[\CL_{A,v,\chi}=\Hom\left(\Gal(L_\chi/F_v),T\right)=\Hom\left(\Gal(E_\chi/F_v),T\right).\]
Since $T$ is annihilated  by $p$, the homomorphisms factor through $\Gal(L_\chi/F_v)\otimes \BZ/p\BZ=\Gal(E_\chi/F_v)$, and the last equality follows as desired.
\end{proof}

Suppose $\omega\in \BO_\bp$, $v\in \CP_1$ and $\chi\in \CC_{n}^\tr(F_v)$. Fix an $\BO_\bp$-generator $x_{d+1}\in A[\bp^{d+1}]$ and define a  map $\xi_\chi: G_v\ra T$ by, for all $\sigma\in G_v$,
\[ \xi_\chi(\sigma)=(\chi\theta_v(\sigma)-1)x_{d+1}.\]

\begin{prop}\label{local-condition2}
Suppose $\omega\in \BO_\bp$, $v\in \CP_1$ and $\chi, \mu\in \CC_{n}^\tr(F_v)$.
\begin{itemize}
\item[(1)] The map $\xi_\chi$ defines a $1$-cocyle in $\RH^1(F_v,T)$.
\item[(2)] The $1$-cocycle $\xi_\chi$ is ramified and generates $\CL_{A,v,\chi}$.
\item[(3)]  If $\xi_\chi$ and $\xi_\mu$ generate the same line in $\RH^1(F_v,T)$, then $\xi_\chi=a\xi_\mu$ for some $a\in \BF_p^\times$, or equivalently $\mu^a\equiv \chi \theta^{1-a}_v\mod \bp^{d+1}$.
\item[(4)] $\xi_\chi=\xi_\mu$ if and only if $\chi/\mu\in \CC_{\leq n-1}(F_v)$.
\end{itemize}
\end{prop}
\begin{proof}
By Lemma \ref{gen}.(2),  $1\otimes x_{d+1}$ is a generator of $A_\chi[\fp^{d+1}]$ and $A_\chi[\fp^{d+1}](F_v)=A_\chi[\fp^d]$. Since $v\nmid p\infty$, $A_\chi(F_v)=A_\chi[\fp^d]\times B$ where $B$ is a profinite group such that $\fp B=B$. Then the image  of the local Kummer descent map
\[\kappa_{A_\chi, \fp,v}: A_\chi(F_v)/\fp A_\chi(F_v)\hookrightarrow \RH^1(F_v, A_\chi[\fp])\]
is generated by the cocycle  $c\in \RH^1(F_v, A_\chi[\fp])$ where, for $\sigma\in G_v$,
\begin{equation}\label{transfer}
c(\sigma)=(1\otimes x_{d+1})^{\sigma-1}=\chi(\sigma)\otimes \theta_v(\sigma)x_{d+1}-1\otimes x_{d+1}=1\otimes (\chi\theta_v(\sigma)-1)x_{d+1}.
\end{equation}
Thus under the identification $A_\chi[\fp]\cong A[\bp]$ of Theorem \ref{twist-av}, the cocycle $c$ transfers to $\xi_\chi\in \RH^1(F_v, T)$ and $\CL_{A,v,\chi}$ is generated by $\xi_\chi$. Since $A$ has good reduction at $v$, $\theta_v$ is unramified, while $\chi$ is totally ramified of order $p^n$. As the proof of Lemma \ref{gen}.(2), one proves $\xi_\chi$ is ramified. This completes the proof of (1) and (2). 

Since $v\in \CP_1$, $G_v$ acts trivially on $T$. Under local class field theory, one has
\[\RH^1(F_v,T)=\Hom(G_v,T)=\Hom(F_v^\times/F_v^{\times  p}, T).\]
Since  $v\nmid p\infty$, 
\[F_v^\times/F_v^{\times p}=\langle \varpi_v \rangle^{\BZ/p\BZ}\times \langle t\rangle^{\BZ/p\BZ},  \]
where $\varpi_v$ is a uniformizer of $F_v$ and $t$ is a generator of $\CO_{v}^\times/\CO_v^{\times p}$. Thus for $\chi\in \CC_{n}^\tr(F_v)$, $\xi_\chi$ is uniquely determined by its values at $\varpi_v$ and $t$.

For $\chi,\mu\in \CC_{n}^\tr(F_v)$, suppose $\xi_\chi$ and $\xi_\mu$ generate the same line, i.e. $\xi_\chi=a \xi_\mu$ for some $a\in k=\BO/\bp$. Equivalently,  for all $s\in F_v^\times$,
\begin{equation}\label{mod-iden}
(\chi(s)\theta_v(s)-1)\equiv a (\mu(s)\theta_v(s)-1)\mod \bp^{d+1}.
\end{equation}
By the definition of $d$, $\omega$ and hence $\chi(s),\mu(s)\in 1+\bp^{d}\BO_\bp$. Since $v\in \CP_1$, $G_v$ acts trivially on $A[\bp^d]$, and thus $\theta_v(s)\in 1+\bp^d\BO_\bp$. Through the natural group isomorphism
\begin{equation}\label{mod-isom}
(1+\bp^d\BO_\bp)/(1+\bp^{d+1}\BO_\bp)\ra \bp^d/\bp^{d+1}, \quad 1+x\mapsto x\mod \bp^{d+1},
\end{equation}
we conclude if $a\in \BF_p^\times$, then
\[(\chi(s)\theta_v(s)-1)\equiv a (\mu(s)\theta_v(s)-1)\equiv ((\mu(s)\theta_v(s))^a-1)\mod \bp^{d+1},\]
and equivalently $\chi\theta_v\equiv (\mu\theta_v)^a\mod \bp^{d+1}$.

Next we prove $a\in \BF_p^\times$. Since both $\chi$ and $\mu$ are totally ramified of order $p^n$, both $\chi(t)$ and $\mu(t)$ are  primitive $p^n$-th roots of unity by Lemma \ref{character}. Then there exists $i\in (\BZ/p^n\BZ)^\times$ such that $\chi(t)=\mu(t)^i$.  Since $\theta_v$ is unramified, $\theta_v(t)=1$.  By (\ref{mod-iden}) and (\ref{mod-isom}),
\begin{equation*}\label{equality1}
(\chi(t)-1)\equiv (\mu(t)^i-1)\equiv i(\mu(t)-1)\equiv a(\mu(t)-1)\mod \bp^{d+1}.
\end{equation*}
It follows that $(a-i)(\mu(t)-1)\in \bp^{d+1}\BO_\bp$. Since $\mu(t)-1\notin \bp^{d+1}\BO_\bp$, necessarily we have $a\equiv i\mod \bp$, i.e. $a\in \BF_p^\times$. This proves (3), and (4) follows immediately.
\end{proof}

Through the Kummer $\fp$-descents, we have a map
\[\CL_{A,v}: \CC(F_v)\lra \CI_v, \quad \chi\mapsto \CL_{A,v,\chi}.\]
\begin{prop}\label{equi}
\begin{itemize}
\item[(1)] For $v\in \CP_0$, the map $\CL_{A,v}$ is constant with $\RH^1_\ur(F_v,T)$ as its image.
\item[(2)] Suppose $v\in \CP_1$. The restriction $\CL_{A,v}:\CC_n^\tr(F_v)\ra \CI_v^\ram$ is surjective and each fibre has the same cardinality $p^{2n-2}(p-1)$.
\end{itemize}
\end{prop}
\begin{proof}
Assertion (1) follows from Proposition \ref{local-condition1}. Suppose $v\in \CP_1$. By Proposition \ref{local-condition3} and \ref{local-condition2}, for $\chi\in C_n^\tr(F_v)$, $\CL_{A,v,\chi}$ is ramified, and each fibre of $\CL_{A,v}:\CC_n^\tr(F_v)\ra \CI_v^\ram$ has cardinality (at most) $|C_{\leq n-1}(F_v)|\cdot |\BF_p^\times|=p^{2n-2}(p-1)$. Since 
\[|\CC_n^\tr(F_v)|=p^{2n-1}(p-1),\quad |\CI_v^\ram|=p,\]
we conclude that $\CL_{A,v}:\CC_n^\tr(F_v)\ra \CI_v^\ram$ must be surjective with fibers of equal cardinality $p^{2n-2}(p-1)$. This finishes the proof.
\end{proof}

\subsection{Twisted incoherent Selmer structures}
\begin{defn}\label{Selmer group}
\begin{itemize}
\item[(1)] A Selmer structure $\CS=\{\CL_{\CS,v}\}_v$ is a collection of subsapces $\CL_{\CS,v}\subset \RH^1(F_v,T)$ for each $v$, such that for almost all $v$, $\CL_{\CS,v}=\RH^1_\ur(F_v,T)$. 
\item[(2)] For a Selmer structure $\CS$, define the associated Selmer group as
\[\RH^1_\CS(F,T)=\Ker\left(\RH^1(F,T)\ra\prod_v\RH^1(F_v,T)/\CL_{\CS,v}\right).\]
\end{itemize}
\end{defn}

Let $\chi\in \CC(F)$ be a global character and for each $v\in \Sigma$ denote by $\chi_v$ the restriction of $\chi$ to $G_v$. Through the local Kummer descents on the abelian variety $A_\chi$, one has a Selmer structure 
\[\{\CL_{A,v,\chi_v}: \text{$v$ runs through all places of $F$}\}\]
whose associated Selmer group is the $\fp$-Selmer group $\Sel_{\fp}(A_\chi).$

Recall $\CB$ is a finite set of places of $F$ containing all the archimedean places, all the places above $p$ and all the places where $A$ has bad reductions, and $\CP_1=\{v\notin\CB: F_v(A[\bp^d])=F_v\}.$
Define
\[\CD=\{\text{squarefree products of primes }v\in  \CP_1\}.\]
For $\delta\in \CD$, define
\[\Omega_\delta:=\prod_{v\in \CB} \CC(F_v)\times \prod_{v\mid \delta} \CC_{n}^\tr(F_v).\]
For $v\in \CP_1\backslash \delta$, by which we mean that $v\in \CP_1$ but  $v\nmid \delta$, 
let $\eta_{\delta,v}:\Omega_{\delta v}\ra \Omega_\delta$ be the natural projection map. Denote $\CB(\delta)=\CB\bigcup \{v:v\mid \delta\}$.

\begin{defn}\label{incoherent Selmer}
For $\omega=(\omega_v)\in \Omega_\delta$, define the twisted  incoherent  Selmer structure $\CS(\omega)$ so that $\CL_{\CS(\omega),v}=\CL_{A,v,\omega_v}$ for $v\in \CB(\delta)$ and $\CL_{\CS(\omega),v}=\RH^1_\ur(F_v,T)$ for $v\notin \CB(\delta)$.
\end{defn}
For $\delta\in \CD$ and $\omega\in \Omega_\delta$, denote
\[\Sel(T,\omega):=\RH^1_{\CS(\omega)}(F,T),\quad r(\omega):=\dim_{k} \Sel(T,\omega),\]
and 
\[t(\omega,v):=\dim_{k} \loc_v\left(\Sel(T,\omega)\right),\]
where $\loc_v: \RH^1(F,T)\ra \RH^1(F_v,T)$ is the localization map at $v$.

\begin{prop}\label{selmer-change}
Suppose $\delta\in \CD$, $\omega\in \Omega_\delta$ and $v\in \CP_1\backslash \delta$ and $\omega'\in  \eta_{\delta, v}^{-1}(\omega)$. Then $0\leq t(\omega,v) \leq 1$, and 
\[r(\omega')-r(\omega)=\left\{\begin{aligned}
+1,\quad & \text{if $t(\omega,v)=0$}, \text{ for exactly $p^{2n-2}(p-1)$ characters $\omega'$};\\
0,\quad & \text{if $t(\omega,v)=0$}, \text{ for exactly $p^{2n-2}(p-1)^2$ characters $\omega'$};\\
-1,\quad & \text{if $t(\omega,v)=1$}.
\end{aligned}\right.\]
\end{prop}
\begin{proof}
Define the relaxed and strict Selmer groups associated to $\CS(\omega)$ at $v$ respectively as
\[\Sel(T,\omega)^{(v)}:=\Ker\left(\RH^1(K,T)\xrightarrow{\sum_{w\neq v} \loc_w}\bigoplus_{w\neq v}\RH^1(K_w,T)/L_{\CS(\omega),w}\right),\]
\[\Sel(T,\omega)_{(v)}:=\Ker\left(\Sel(T,\omega)^{(v)}\xrightarrow{\loc_v}\RH^1(F_v,T)\right).\]
Let $V:=\loc_v(\Sel(T,\omega)^{(v)})\subset \RH^1(F_v,T)$. By the Poitou-Tate global duality (see for example  \cite[Theorem 3.1]{Tateduality}, \cite[Theorem I-4.10]{Milne-ADT06} and \cite[Theorem 1.10.2]{Rubin}), $V$ is maximal isotropic in $\RH^1(F_v,T)$, and hence $\dim_{k}V=1$. Since $\loc_v(\Sel(T,\omega))=V\cap \RH^1_\ur(F_v,T)$, $0\leq t(\omega,v)\leq 1$.

We have the following exact sequences
\[0\ra \Sel(T,\omega)_{(v)}\ra \Sel(T,\omega)\ra V\cap \RH^1_\ur(F_v,T)\ra 0,\]
\[0\ra \Sel(T,\omega)_{(v)}\ra \Sel(T,\omega')\ra V\cap \CL_{\CS(\omega'),v}\ra 0.\]
Then
\[r(\omega')-r(\omega)=\dim_{k}(V\cap \CL_{\CS(\omega'),v})-t(\omega,v).\]

If $t(\omega,v)=1$ then $V=\RH^1_\ur(F_v,T)$. Since $\omega'_v$ is totally ramified of order $p^n$, by Proposition \ref{local-condition2}, $\CL_{\CS(\omega'),v}=\CL_{A,v,\omega'_v}$ is ramified. Then   $V\cap \CL_{\CS(\omega'),v}=0$ and $r(\omega')-r(\omega)=-1$. 

Suppose $t(\omega,v)=0$. Then $V\cap \RH^1_\ur(F_v,T)=0$. If  $V= \CL_{\CS(\omega'),v}$, then $r(\omega')-r(\omega)=1$; Otherwise $r(\omega')-r(\omega)=0$. The proposition follows by counting the number of local characters $\omega'_v\in \CC_{n}^\tr(F_v)$ whose local Kummer images $\CL_{A,v,\omega'_v}$ coincide with $V$ according to Proposition  \ref{equi}.
\end{proof}

For $\delta\in \CD$, define 
\[\Delta(\delta):=\left\{
\chi\in \CC(F):\begin{aligned} &\text{$\chi$ is totally ramified of order $p^n$ 
at all $v\mid \delta$ }\\&\text{ and unramified outside $\CB(\delta)\bigcup \CP_0$}\end{aligned}\right\}.\]
\begin{prop}\label{Selmer}
For $\delta\in \CD$ and a character $\chi\in \Delta(\delta)$, one has 
\[\Sel_{\fp}(A_\chi)=\Sel(T,\omega)\] 
where $\omega=(\chi_v)_{v\in \CB(\delta)}\in \Omega_\delta$.
\end{prop}
\begin{proof}
By definition,
\[\Sel_{\fp}(A_\chi)=\{x\in \RH^1(F,T): x\in \CL_{A,v,\chi_v} \text{ for all }v\}.\]
We check the local conditions for each $v$. For $v\in \CB(\delta)$, by definition, $\CL_{\CS(\omega),v}=\CL_{A,v,\chi_v}$. For $v\notin \CB(\delta)$, either $v\in \CP_0$ or that $v\in \CP_1$ and $\chi_v$ is unramified. In both cases, according to Proposition \ref{local-condition1},
\[\CL_{A,v,\chi_v}=\RH^1_\ur(F_v,T)=\CL_{\CS(\omega),v}\]
as desired.
\end{proof}

\section{Markov model for twisted incoherent Selmer groups}\label{Markov-process}
In this section we prove that the rank distribution of the twisted incoherent Selmer groups is governed by certain Markov process using the machinery of Klagsbrun-Mazur-Rubin \cite[Part I]{KMR14}. We first record the following effective Chebotarev theorem that we will use later and which can be found in \cite[Theorem 2 and 4]{Serre81} and \cite[Theorem 8.1]{KMR14}.
\begin{thm}\label{chebotarev}
There is a nondecreasing continuous function $\CL:[1,\infty]\ra [1,\infty]$ such that for 
\begin{itemize}
\item every $Y\geq 1$,
\item every $\delta\in \CD$ with $\N_\delta<Y$, where $\N_\delta$ denotes the norm of $\delta$,
\item every Galois extension $L$ of $F$ that is abelian of exponent $p^n$ over $F_1=F(T)$, and unramified outside $\CB(\delta)$,
\item every pair of subsets $S,S'\subset \Gal(L/F)$ stable under conjugation, with  $S$ nonempty, and 
\item every $X>\CL(Y)$,
we have
\[\left|\frac{|\{v\notin\CB(\delta): \N_v\leq X, \Fr_v(L/F)\in S'\}|}{|\{v\notin\CB(\delta): \N_v\leq X, \Fr_v(L/F)\in S\}|}-\frac{|S'|}{|S|}\right|\leq \frac{1}{Y}.\]
In particular, $\{v\notin\CB(\delta): \N_v\leq X, \Fr_v(L/F)\subset S\}$ is nonempty.
\end{itemize}
\end{thm}

As a subgroup of $\Aut_{k}(T)\cong k^\times$, $G=\Gal(F_1/F)$ has order prime to $p$. Then for any $i\geq 0$, $\RH^i(G,T)=0$, thus the inflation-restriction exact sequence gives an isomorphism 
\[\Res_{F_1/F}: \RH^1(F,T)\xrightarrow{\cong} \Hom_G(G_{F_1},T).\]
For $\delta\in\CD$ and $\omega\in \Omega_\delta$, denote $W=\Res_{F_1/F}(\Sel(T,\omega))$. Therefore we have a $G$-equivariant pairing
\begin{equation}
\label{pairing1} G_{F_1}\times W\ra T.
\end{equation}
Let $F_\omega$ be the fixed field of  the left kernel of the above paring. Since cocycles in $\Sel(T,\omega)$ are unramified outside $\CB(\delta)$, $F_\omega/F_1$ is unramifed outside $\CB(\delta)$. 

\begin{prop}\label{control-field}
We have
\[\Gal({F_{\omega}}/F_1)\cong \Hom_k(\Sel(T,\omega),T)\text{ and }\Sel(T,\omega)\cong \Hom_{\BF_p[G]}(\Gal({F_{\omega}}/F_1),T).\]
\end{prop}
\begin{proof}
Under the pairing (\ref{pairing1}), we have an injection \[\Gal({F_{\omega}}/F_1)\hookrightarrow \Hom_k(W,T)\cong T^{r(\omega)}\] as $\BF_p[G]$-modules. By Hypothesis ({\bf H2}), $T$ is a simple $\BF_p[G]$-module. As an $\BF_p[G]$-submodule of the semisimple module $T^{r(\omega)}$, $\Gal({F_{\omega}}/F_1)\cong T^i$ for some $0\leq i\leq r(\omega)$.  Since $\Hom_{\BF_p[G]}(T,T)=\Hom_k(T,T)=k$, 
\[\Hom_{\BF_p[G]}(\Gal({F_{\omega}}/F_1),T)\cong\Hom_k(T^i,T)=k^i.\]
Through the pairing (\ref{pairing1}), $W$ injects into $\Hom_{\BF_p[G]}(\Gal({F_{\omega}}/F_1),T)$, and  thus has $k$-dimension $\leq i$. Then we conclude $i=r(\omega)$, and the proposition follows as desired.
\end{proof}

\begin{defn}\label{constant}
For $\gamma\in \Omega_1$, set
\[\Omega^\gamma_\delta:=\{\gamma\}\times \prod_{v\mid \delta} \CC_{n}^\tr(F_v),\]
and 
\[\Sha_{\CP_1}(\gamma)=\{s\in \Sel(T,\gamma): s_v=0 \text{ for all }v\in \CP_1\},\quad r_\gamma=\dim_{k} \Sha_{\CP_1}(\gamma).\]
\end{defn}
Clearly $\Omega_\delta=\bigsqcup_{\gamma\in \Omega_1}\Omega_\delta^\gamma$, and 
 for any $\omega\in \Omega^\gamma_\delta$, 
\begin{equation*}\label{lower-bound}
\Sha_{\CP_1}(\gamma)\subset \Sel(T,\omega),\quad r_\gamma\leq r(\omega).
\end{equation*}
Let $F_{\Sha_{\CP_1}(\gamma)}$ be the fixed field of the annihilator of $\Sha_{\CP_1}(\gamma)$ under the paring (\ref{pairing1}).
\begin{prop}\label{disjoint}
 For any $\omega\in \Omega^\gamma_\delta$, we have $F_\omega\cap F_d=F_{\Sha_{\CP_1}(\gamma)}$ and the isomorphisms in Proposition \ref{control-field} induce an isomorphism $$\Gal(F_{\Sha_{\CP_1}(\gamma)}/F_1)\cong \Hom_k(\Sha_{\CP_1}(\gamma),T).$$
\end{prop}
\begin{proof}
The isomorphism follows immediately from Proposition \ref{control-field}. For $F_{\Sha_{\CP_1}(\gamma)}=F_\omega\cap F_d$, it suffices to prove that for $s\in \Sel(T,\omega)$, $s\in \Sha_{\CP_1}(\gamma)$ if and only if $s(G_{F_d})=0$. For any $v\in \CP_1$, $v$ is split in $F_d/F$, i.e. $ G_v \subset G_{F_d}$. If  $s(G_{F_d})=0$, then  $s_v=0$. 

Next suppose $s\in \Sha_{\CP_1}(\gamma)$ and $\sigma\in G_{F_d}$. By Chebotarev density theorem, there exists a place $v\notin \CB(\delta)$ such that $\Fr_v=\sigma$. Then $\Fr_v$ acts trivially on $F_d$, i.e. $v$ belongs to $\CP_1$. By the definition of $\Sha_{\CP_1}(\gamma)$, $s(\sigma)=s(\Fr_v)=0$, and hence $s(G_{F_d})=0$.
\end{proof}

For $\omega\in \Omega^\gamma_\delta$, to justify the presence of $\Sha_{\CP_1}(\gamma)$ in $\Sel(T,\omega)$, we define the normalized rank map 
\[r^\gamma: \Omega_\delta^\gamma\ra \BZ_{\geq 0},\quad r^\gamma(\omega)=r(\omega)-r_\gamma.\]
For $i=0,1$, denote
\[c_i(\omega)=\left\{\begin{aligned}
&q^{-r^\gamma(\omega)},\quad && \text{if $i=0$},\\
&1-q^{-r^\gamma(\omega)},\quad && \text{if $i=1$},
\end{aligned}\right.\]
where $q$ is the cardinality of the residue field $k$.
For $X>0$ define
\[\CP_1(X):=\{v\in \CP_1: \N_v\leq X\} \text{ and }\]
\[\Phi_i (\omega,X):=\frac{|\{v\in \CP_1(X)\backslash \delta: t(\omega,v)=i\}|}{| \CP_1(X)\backslash \delta|}.\]

\begin{prop}\label{probabilityfort}
Suppose $\gamma\in \Omega_1$, $\delta\in \CD$, $\omega\in \Omega^\gamma_\delta$, $v\in \CP_1\backslash \delta$ and $\CL$ is the function as in Theorem \ref{chebotarev}.
For every $Y>\N_\delta$, and every $X>\CL(Y)$
\[\left|\Phi_i(\omega,X)-c_i(\omega)\right|\leq \frac{1}{Y}.\]
\end{prop}
\begin{proof}
For $v\in \CP_1\backslash \delta$, $\Fr_v$ acts trivially on $T$. By \cite[XIII.1]{Serre-localfields}, 
$$\RH^1_\ur(F_v,T)\cong T/(\Fr_v-1)T=T,\quad s\mapsto s(\Fr_v).$$ Thus 
\[t(\omega,v)=\dim_{k} \left\{s(\Fr_v): s\in \Sel(T,\omega)\right\}.\]
Then $t(\omega,v)=0$ if and only if $\Fr_v \in G_{F_\omega}$ by the definition of $F_\omega$. Let $F_\omega F_d$ be the composition field. A place $v\notin \CB(\delta)$ belongs to $\CP_1$ if and only if $\Fr_v|_{F_d}=1$. Put
\[S=\{\sigma\in\Gal(F_\omega F_d/F_{\Sha_{\CP_1}(\gamma)}): \sigma|_{F_d}=1 \},\]
and 
\[S_0=\{1\}, \quad S_1 =S\backslash \{1\}.\]
Then for $v\in \CP_1\backslash \delta$, $t(\omega,v)=i$ if and only if $\Fr_v\in S_i$. By Proposition \ref{disjoint}, $F_\omega\cap F_d=F_{\Sha_{\CP_1}(\gamma)}$ and thus
\[\Gal(F_\omega F_d/F_{\Sha_{\CP_1}(\gamma)})\cong \Gal(F_\omega /F_{\Sha_{\CP_1}(\gamma)})\times \Gal(F_d/F_{\Sha_{\CP_1}(\gamma)}).\]
Then, combined with Proposition \ref{control-field},  \ref{disjoint}, we have
\[ |S|=[F_\omega:F_{\Sha_{\CP_1}(\gamma)}]=q^{r(\omega)-r_\gamma}, \ |S_0|=1\text{ and } |S_1|=q^{r(\omega)-r_\gamma}-1.\]
Applying Theorem \ref{chebotarev}, for $X>\CL(Y)$ we have
\[\left|\frac{|\{v\in \CP_1(X)\backslash \delta: t(\omega,v)=i\}|}{| \CP_1(X)\backslash \delta|}-\frac{|S_i|}{|S|}\right|\leq \frac{1}{Y}\]
as desired.
\end{proof}

In the rest of this section, we follow the terminologies for Markov process as in \cite[Part I]{KMR14}. Denote $\ell^1$ the Banach space over $\BR$
\[\ell^1:=\{\mu: \BZ_{\geq 0}\ra \BR: \|\mu\|:=\sum_{r\geq 0} |\mu(r)| <\infty\}.\]
\begin{defn}
Define the space $W\subset \ell^1$ of probability distributions to be
\[W:=\{\mu\in \ell^1: \mu(r)\geq 0 \text{ for all $r\in \BZ_{\geq 0}$ and } ||\mu||=1\}.\]
A bounded linear operator $M:\ell^1\ra\ell^1$ is called a Markov operator if $M(W)\subset W$.
\end{defn}
We can view $M$ as an infinite matrix $[m_{r,s}], r,s\in \BZ_{\geq0}$ where, for $\mu\in \ell^1$,
\[(\mu M )(s)= \sum_{r\geq 0} m_{r,s}\mu(r)\]
with $m_{r,s}$ bounded, and then $M$ is a Markov operator if and only if $m_{r,s}\geq 0$ for
all $r,s\geq 0$ and $\sum_{s\geq 0}m_{r,s} = 1$ for every $r$.

We fix  $\gamma\in \Omega_1$. For $\delta\in \CD$, $v\in \CP_1\backslash \delta$, the finite sets  $\Omega_\delta^\gamma$, together with the projections $\eta_{\delta,v}:\Omega^\gamma_{\delta v}\ra \Omega^\gamma_\delta$ and  the normalized rank maps $ r^\gamma:\Omega^\gamma_\delta\ra \BZ_{\geq 0}$ form a {\em rank data} $\Omega^\gamma$ on $\CD$ (see \cite[Definition 3.7]{KMR14}). Here it suffices to note that each fiber $\eta_{\delta,v}^{-1}(\omega)$ has the same cardinality $|\CC_{n}^\tr(F_v)|=p^{2n-1}(p-1)$ which is independent of $\delta\in\CD$, $v\in \CP_1\backslash \delta$ and $\omega\in \Omega^\gamma_\delta$. Then we have the rank distribution function
\[E^\gamma_\delta(r):=\frac{|\{\omega\in \Omega^\gamma_\delta: r^\gamma(\omega)=r\}|}{|\Omega^\gamma_\delta|},\quad r\in \BZ_{\geq 0}.\]
If $B$ is a nonempty finite subset of $\CD$, the rank distribution over $B$
is the average of the $E^\gamma_\delta$ over $\delta\in B$, weighted according to the size of 
$\Omega^\gamma_\delta$
\[E_B^\gamma:=\frac{\sum_{\delta\in B}  |\Omega_\delta^\gamma|E_\delta^\gamma}{\sum_{\delta\in B} |\Omega^\gamma_\delta|}.\]
Thus $E^\gamma_B(r)$ is the probability of, as $\delta$ ranges through $B$, that $r^\gamma(\omega)=r$. If  for all $\delta\in B$, $\Omega^\gamma_\delta$ have the same cardinality, then $E^\gamma_B=\frac{\sum_{\delta\in B}E_\delta^\gamma}{|B|}$.

Denote
\[\epsilon=\left\{\begin{aligned}
&\sym=1,\quad && \text{if $(\CMS)$ holds},\\
&\uni=1/2,\quad && \text{if $(\CMU)$ holds}.
\end{aligned}\right.\]
If $(\CMS)$ resp. $(\CMU)$ holds, the residue field $k=\CO/\fp=\BF_q$ with $q=p$ resp. $q={p^2}$. Define a  Markov operator $M_\epsilon=[m^\epsilon_{r,s}]$ on $\ell^1$ by 
\begin{equation}\label{markov}
m^\epsilon_{r,s}=\left\{\begin{aligned}
&1-q^{-r},\quad && \text{if $s=r-1\geq 0$},\\
&(1-q^{-\epsilon})q^{-r},\quad && \text{if $s=r\geq 0$},\\
&q^{-r-\epsilon},\quad && \text{if $s=r+1\geq 1$},\\
&0,\quad && \text{otherwsie}.
\end{aligned}\right.
\end{equation}

\begin{thm}
For $\gamma\in \Omega_1$, the rank data $\Omega^\gamma=(\Omega_\delta^\gamma, \eta_{\delta,v},r^\gamma)$ on $\CD$ is governed by the Markov operator $M_\epsilon$ with convergence rate $\CL$ in the sense of \cite[Definition 3.10-12]{KMR14}.
\end{thm}
\begin{proof}
Let $\omega\in \Omega_\delta^\gamma$ and $r^\gamma(\omega)=r$.
Define
\[F(\omega,X,s):=\frac{\sum_{{v\in \CP_1(X)\backslash \delta}}|\{\omega'\in \eta_{\delta,v}^{-1}(\omega): r^\gamma(\omega')=s\}|}{\sum_{{v\in \CP_1(X)\backslash \delta}}|\eta_{\delta,v}^{-1}(\omega)|},\]
By Proposition \ref{selmer-change},
\[F(\omega,X,s)=\left\{\begin{aligned}
&\Phi_1(\omega,X),\quad && \text{if $s=r-1\geq 0$},\\
&\frac{p-1}{p}\Phi_0(\omega,X),\quad && \text{if $s=r\geq 0$},\\
&\frac{1}{p}\Phi_0(\omega,X),\quad && \text{if $s=r+1\geq 1$},\\
&0,\quad && \text{otherwise}.
\end{aligned}\right.\]
Then Theorem \ref{probabilityfort} shows 
\[\lim_{X\ra \infty} F(\omega,X,s)=\left\{\begin{aligned}
&1-q^{-r},\quad && \text{if $s=r-1$},\\
&(1-q^{-\epsilon})q^{-r},\quad && \text{if $s=r$},\\
&q^{-r-\epsilon},\quad && \text{if $s=r+1$},\\
&0,\quad && \text{otherwise}.
\end{aligned}\right.\]
The precise convergence in Theorem \ref{probabilityfort} shows that $\CL$ is a convergence rate for $(\Omega^\gamma,M_\epsilon)$.

\end{proof}

Define the symplectic ($\epsilon=1$) or unitary ($\epsilon=1/2$) distribution  $\sD_q^\epsilon:\BZ_{\geq 0}\ra \BR$ as
\[\sD_q^\epsilon(r):=\prod_{i\geq 0}(1+q^{-i-\epsilon})^{-1}\prod_{i=1}^r\frac{q^{1-\epsilon}}{q^i-1}.\]

\begin{prop}\label{distribution}
Suppose $\epsilon=1$ or $1/2$.
\begin{itemize}
\item[(1)]  We have 
\[\sum_{r\geq 0} r\sD_q^\epsilon(r)=\sum_{i=0}^\infty \frac{1}{1+q^{i+\epsilon}},\]
\[\sum_{r\geq 0} q^r\sD_q^\epsilon(r)=1+q^{1-\epsilon},\]
and 
\[\sum_{r\text{ odd}} \sD^\epsilon_q(r)=\frac{1-\beta^\epsilon}{2}\text{ where }\beta^\epsilon=\prod_{i=0}^\infty\frac{1-q^{-i-\epsilon}}{1+q^{-i-\epsilon}}.\]
\item[(2)] We have 
\[\sD_q^\epsilon(0)=\prod_{i\geq 0}(1+q^{-i-\epsilon})^{-1}>1-\frac{q^{1-\epsilon}}{q-1}.\]
In particular, $\lim\limits_{q\ra \infty}\sD_q^\epsilon(0)=1$.
\end{itemize}
\end{prop}
\begin{proof}
See \cite[Proposition 2.21]{Shu-quadratic}. 
\end{proof}

\begin{prop}\label{lambda}
The Markov operator $M_\epsilon$ in $(\ref{markov})$ has invariant distribution $\sD_q^\epsilon$, and for every distribution $\mu\in W$, one has
\[\lim_{k\ra \infty} \mu M_\epsilon^k=\sD_q^\epsilon.\]
\end{prop}
\begin{proof}
It is straightforward to verify, for both the $(\CMS)$ and $(\CMU)$ cases, that $\sD_q^\epsilon$ is invariant for the Markov operator $M_\epsilon$, i.e. $\sD_q^\epsilon M_\epsilon=\sD_q^\epsilon$. We view $M_\epsilon$ as a Markov  process on $\BZ_{\geq 0}$. This Markov process is irreducible and aperiodic in the sense of \cite[Chapter 1]{Norris97}. By \cite[Theorem 1.8.3]{Norris97}, $\sD_q^\epsilon$ is the unique invariant distribution for $M_\epsilon$ and for every distribution $\mu\in W$, one has
\[\lim_{k\ra \infty} \mu M_\epsilon^k=\sD_q^\epsilon\]
as desired.
\end{proof}

\begin{defn}\label{fan}$($cf. \cite[Definition 3.14]{KMR14}$)$
Define a sequence of real valued functions $\{L_i(X)\}_{i\geq 1}$ by
\begin{eqnarray*}
&&L_1(X)=\CL(X),\\
&&L_{i+1}(X)=\max\left(\CL\left(\prod_{j\leq i}L_j(X)\right),XL_i(X)\right), \quad i\geq 1.
\end{eqnarray*}
For $k\in \BZ_{\geq 0}$ and $X>0$, define 
\[\CD_k(X)=\{\delta\in \CD:\delta=v_1v_2\cdots v_k \text{ with } \N_{v_i}<L_i(X) \text{ for all }i\},\]
and $\CD_k=\bigcup_{X>0} \CD_k(X)$. 
\end{defn}

\begin{thm}\label{local-markov}
Suppose $\gamma\in \Omega_1$. For every $k\geq 0$, 
\[\lim_{X\ra \infty} E^\gamma_{\CD_k(X)}=E^\gamma_1M_\epsilon^k,\]
and 
\[\lim_{k\ra \infty}\lim_{X\ra \infty} E^\gamma_{\CD_k(X)}=\sD_q^\epsilon.\]
\end{thm}
\begin{proof}
By Proposition \ref{selmer-change}, for any $k\geq 0$, $\delta\in \CD_k$ and $\omega\in \Omega^\gamma_\delta$, $r(\omega)\leq k+r_\gamma$. Thus the hypothesis of  \cite[Theorem 4.3]{KMR14} holds. Thus the theorem follows from   \cite[Theorem 4.3]{KMR14}  and Proposition \ref{lambda}.
\end{proof}

\section{Passage to global characters}
In this section, we build a bridge between the global characters and the incoherent systems of local characters so that the statistical behavior of the twisted incoherent Selmer groups established in Section \ref{Markov-process} will be transferred to  Selmer groups twisted by global characters through this bridge in the next section. First we define the subset of global characters we are interested in as follows.

\begin{defn}\label{charset}
Let $\Delta(F)\subset \CC(F)$ be the subset of characters $\chi$ that if $\chi_v$ is ramified at a place $v\in \CP_1$ then $\chi_v$ is totally ramified of order $p^n$.  
\end{defn}
In particular, if $n=1$, then $\Delta(F)=\CC(F)$. For $\delta\in \CD$, define 
\[\Delta(\delta):=\left\{
\chi\in \CC(F):\begin{aligned} &\text{$\chi$ is totally ramified of order $p^n$ 
at all $v\mid \delta$ }\\&\text{ and unramified outside $\CB(\delta)\bigcup \CP_0$}\end{aligned}\right\}.\]
Since $\Omega=\CB\bigsqcup \CP_0\bigsqcup \CP_1$, $\Delta(F)=\bigsqcup_{\delta\in \CD}\Delta(\delta)$.

 In what follows we assume $\Pic(\CO_{\CB})=0$ which can always be achieved by suitably enlarging $\CB$. Then global class field theory gives
\begin{eqnarray}\label{lg}
\begin{split}
\CC(F)&=&&\Hom(F^\times\backslash \BA_F^\times,\mu_{p^n})&\\
&=&&\Hom\left(\left.\left(\prod_{v\in \CB(\delta)}F_v^\times\times \prod_{v\notin \CB(\delta)}\CO_v^\times\right)\right/\CO_{\CB(\delta)}^\times,\mu_{p^n}\right).&
\end{split}
\end{eqnarray}

\subsection{Case $F\neq F_d$: Enough free places}
In this subsection, we assume $F\neq F_d$. Then $\CP_0$ is the set of places $v\notin \CB$ nonsplit in $F_d/F$ and $\CP_0$ has positive density.  By definition we have no extra restrictions for characters in $\Delta(F)$ over the places in $\CP_0$ other than the compatibility from global class field theory.  This freeness allows us to choose suitable local characters over $\CP_0$, under the compatibility from global class field theory, to fill the defect between the global characters and the incoherent systems of local characters.

For any $X>0$, define
\[\Delta(X)=\{\chi\in \Delta(F): \text{$\chi$ is unramified outside $\CB\bigcup \{v: N_v<X\}$}\},\]
and set $\Delta(\delta,X)=\Delta(\delta)\cap \Delta(X)$.  Note $\Delta(\delta,X)$ is finite and is nonempty  if $X>\N_\delta$. Let $\ell_{\delta,X}: \Delta(\delta,X)\ra \Omega_\delta$ be the natural restriction map $\chi\mapsto (\chi_v)_{v\in \CB(\delta)}.$

\begin{lem}\label{free}
Suppose $\CL$ is the function as in Theorem \ref{chebotarev}, $\delta\in \CD$, $\alpha\in \CO_{\CB(\delta)}^\times/\CO_{\CB(\delta)}^{\times p^n}$, and $\alpha\neq 1$. Then there exists $v\in \CP_0$ with $\N_v\leq \CL(\N_\delta)$ such that $\alpha\notin \CO_v^{\times p^n}$.
\end{lem}
\begin{proof}
If $n=1$, by the Galois equivariance of  the Weil pairing, $\mu_p\subset F_1\subset F_d$, and if $n\geq 2$, by our assumption, $\mu_{p^n}\subset F$. Thus the field $F(\sqrt[p^n]{\alpha},A[\bp^d])=F(\mu_{p^n}\sqrt[p^n]{\alpha},A[\bp^d])$ is Galois over $F$ and abelian of exponent $p^n$ over $F_1$. In the following we will construct an element $\sigma\in \Gal(F(\sqrt[p^n]{\alpha},A[\bp^d])/F)$ such that both $\sigma|_{F(\sqrt[p^n]{\alpha})}\neq 1$ and $\sigma|_{F_d}\neq 1$. 

If $\sqrt[p^n]{\alpha}\in F_d$, it suffices to take $\sigma\in \Gal(F_d/F)$ acting nontrivially on $\sqrt[p^n]{\alpha}$. Suppose $\sqrt[p^n]{\alpha}\notin F_d$, and take $\sigma\in \Gal(F(\sqrt[p^n]{\alpha},A[\bp^d])/F)$ which maps nontrivially to both the factors of the isomorphism
\[\Gal(F(\sqrt[p^n]{\alpha},A[\bp^d])/F)\cong\Gal(F(\sqrt[p^n]{\alpha},A[\bp^d])/F_d)\rtimes \Gal(F_d/F).\]

Applying Theorem \ref{chebotarev} with $L=F(\sqrt[p^n]{\alpha},A[\bp^d])$ and $S$ the conjugacy class of $\sigma$, there exists a place $v\notin \CB(\delta)$ with $\N_v\leq \CL(N_\delta)$ whose Frobenius in $\Gal(F(\sqrt[p^n]{\alpha},A[\bp^d])/F)$ is the conjugacy class of $\sigma$. For such a  place $v$, we have $v\in \CP_0$ and $\alpha \notin \CO_v^{\times p^n}$.
\end{proof}

\begin{prop}\label{global-local}
Suppose $\delta\in \CD$, $X>\CL(\N_\delta)$. The map $\ell_{\delta,X}:\Delta(\delta,X)\ra \Omega_\delta$ is surjective, and  for every $\omega\in \Omega_\delta$,
\[\frac{|\ell_{\delta,X}^{-1}(\omega)|}{|\Delta(\delta,X)|}=\frac{1}{|\Omega_\delta|}.\]
\end{prop}

\begin{proof}
Set 
\[\CB_1=\{v\in\CP_1:v\nmid \delta \}\bigcup \{v\in \CP_0:\N_v>X\},\]
\[\CB_2=\{v\in \CP_0: \N_v\leq X\}.\]
Recall $\Sigma$ is the set of all places of $F$, and then $\Sigma=\CB(\delta)\bigsqcup \CB_1\bigsqcup \CB_2$. Noting $X>\CL(\N_\delta)$,  $\chi\in \Delta(\delta,X)$ if and only if 
\[(\chi_v)_{v\in \CB(\delta)}\in \Omega_\delta \text{ and }\chi\left(\CO_v^\times\right)=1 \text{ for all $v\in \CB_1$}.\]
The lemma \cite[Lemma 10.4]{KMR14} follows verbatim with $p$ replaced by $p^n$, and apply it with
\[G=\prod_{v\in \CB_2} \CO_v^\times, \quad H=\prod_{v\in \CB(\delta)}F_v^\times\times \prod_{v\in \CB_1}\CO_v^\times, \quad J=\CO_{\CB(\delta)}^\times.\]
Since $X>\CL(\N_\delta)$, by \cite[Lemma 10.4]{KMR14} and Lemma \ref{free}, via (\ref{lg}), the restriction map
\[\CC(F)\ra \Hom\left(\prod_{v\in \CB(\delta)}F_v^\times\times \prod_{v\in \CB_1}\CO_v^\times,\mu_{p^n}\right)\]
is surjective. Then for any $\omega\in \Omega_\delta$ there exists a character $\chi\in \CC(F)$, unramified outside $\CB(\delta)$ and $\CB_2$, which restricts to $\omega$. Such a $\chi$ necessarily belongs to $\Delta(\delta,X)$, and this proves that $\ell_{\delta,X}$ is surjective.

If $\chi_1,\chi_2\in \Delta(\delta,X)$ and $\ell_{\delta,X}(\chi_1)=\ell_{\delta,X}(\chi_2)$, then $\chi_{1}\chi_2^{-1}\in \ell_{\delta,X}^{-1}(1)$. Since $\Delta(\delta,X)$ is stable under multiplication by $ \ell_{\delta,X}^{-1}(1)$, each fibre of $\ell_{\delta,X}: \Delta(\delta,X)\ra \Omega_\delta$ has the same cardinality $| \ell_{\delta,X}^{-1}(1)|$. This completes the proof.
\end{proof}

\subsection{Case $F=F_d$: No free places}
In this subsection we assume $F=F_d$ whence $\CP_0=\emptyset$. In this case, there is no free places and generally the natural restriction map $\ell_{\delta}:\Delta(\delta)\ra \Omega_\delta$ is not surjective due to the compatibility from global class field theory. The key observation is that, when $F=F_d$, for fixed $\delta\in \CD$, any $w\in \CP_1\backslash \delta$ and $\chi\in \Delta(\delta w)$, one always has the exact rank relation 
\[\dim_{k} \Sel_{\fp}(A_\chi)-1=\dim_{k} \Sel(T, \omega)\]
where $\omega=(\chi_v)_{v\in \CB(\delta)}$ is the natural restriction which forgets the place $w$ (see Proposition \ref{dim-preserving}). Such restrictions forgetting one place in $\CP_1$ has two advantages: Firstly they preserve the statistical behavior of Selmer groups thanks to the above exact rank formula. Secondly, varying the forgotten places $w$, they make enough free room passaging from global characters to  incoherent systems of local characters with ``asymptotical equidistributed" fibers and this is achieved by a limiting process in Proposition \ref{fiber} and \ref{global-local2}.

For $\delta=v_1v_2\cdots v_k\in \CD$, we may always arrange the primes with increasing norms so that the decomposition expression is unique. Denote
\[\ov{\delta}=\prod_{v\in \CP_1,\N_v\leq \N_{v_k}} v\text{ and }\Delta(\delta*\CP_1(X))=\bigsqcup_{w\in \CP_1(X)\backslash \ov{\delta}} \Delta(\delta w).\]
Then 
\[\bigsqcup_{\delta'\in \CD_{k+1}}\Delta(\delta')=\lim_{X\ra \infty} \bigsqcup_{\delta\in \CD_k}\Delta(\delta*\CP_1(X)).\]
For the convenience, we introduce the following convergence rate for effective Chebotarev.
\begin{defn}\label{new-convergence}
Let $\CL':[1,\infty]\ra [1,\infty]$ be a nondecreasing continuous function, growing faster than $\CL$, so that, as in Theorem \ref{chebotarev} for every $Y>\N_\delta$ and every $X>\CL'(Y)$, we have
\begin{equation}\label{extended-convergence}
\left|\frac{|\{v\notin\CB(\ov{\delta}): \N_v\leq X, \Fr_v(L/F)\in S'\}|}{|\{v\notin\CB(\ov{\delta}): \N_v\leq X, \Fr_v(L/F)\in S\}|}-\frac{|S'|}{|S|}\right|\leq \frac{1}{Y}.
\end{equation}
\end{defn}
For every $\delta\in \CD$ and every abelian extension $L/F_1$ of exponent $p^n$ unramified outside $\CB(\delta)$.  With Theorem \ref{chebotarev} applied to $\ov{\delta}$ and $L$, for (\ref{extended-convergence}) to hold, it suffices to take $\CL'(\N_\delta)>\CL(\N_{\ov{\delta}})$.

Let $\ell_{\delta,X}: \Delta(\delta*\CP_1(X))\ra \Omega_\delta$ be the restriction map forgetting the places $w\in \CP_1(X)\backslash\ov{\delta}$. Also denote $L_\delta=F\left(\sqrt[p^n]{\CO_{\CB(\delta)}^\times}\right)$. For $\delta\in \CD_k$, the degree $[L_\delta:F]$ is a constant $d(\CB,k)$ depending on $\CB$ and $k$.
\begin{prop}\label{fiber}
For $Y>\N_\delta$ and $X>\CL'(Y)$, the map $\ell_{\delta,X}: \Delta(\delta*\CP_1(X))\ra \Omega_\delta$ is surjective and for any $\omega\in \Omega_\delta$, 
\[\left|\frac{|\ell_{\delta, X}^{-1}(\omega)|}{|\CP_1(X)\backslash \ov{\delta}|}-\frac{p^{n}-p^{n-1}}{d(\CB,k)}\right|\leq \frac{p^n}{Y}.\]
\end{prop}
\begin{proof}
Let $w$ be an arbitrary place in $\CP_1\backslash\ov{\delta}$. Under the isomorphism (\ref{lg}), the characters unramified outside $\CB(\delta w)$ are identified with 
\[\Hom\left(\left.\left(\prod_{v\in \CB(\delta)}F_v^\times \times \CO_w^\times\right)\right/\CO_{\CB(\delta)}^\times,\mu_{p^n}\right).\] 
The restriction map
\begin{equation}\label{lg2}
\Hom\left(\left.\left(\prod_{v\in \CB(\delta)}F_v^\times \times \CO_w^\times\right)\right/\CO_{\CB(\delta)}^\times,\mu_{p^n}\right)\ra \Hom\left(\left.\left(\prod_{v\in \CB(\delta)}F_v^\times\right)\right/U_w,\mu_{p^n}\right),
\end{equation}
is surjective, where $U_w=\Ker\left(\CO_{\CB(\delta)}^\times\ra \CO_w^\times/\CO_w^{\times p^n}\right)$.

For $\omega\in \Omega_\delta$, suppose $\omega\left(\CO_{\CB(\delta)}^\times\right)=\mu_{p^m}$ with $0\leq m\leq n$, and let $U_\omega=\Ker\left(\omega:\CO_{\CB(\delta)}^\times \ra \mu_{p^n}\right)$. If we denote $ L_\omega=F\left(\sqrt[p^n]{U_\omega}\right)$, then $L_\delta/L_\omega$ is cyclic of order $p^m$. Let $\chi$ be a character unramified outside $\CB(\delta)$ which maps to $\omega$ via the restriction map (\ref{lg2}). Then $\chi_w \omega(\CO_{\CB(\delta)}^\times)=1$ and $U_\omega=\CO_{\CB(\delta)}^\times\cap \Ker(\chi_w)$. Thus $\chi\in \Delta(\delta w)$, i.e. $\chi_w$ is totally ramified, if and only if $U_w=U_\omega$. Equivalently, $w$ splits in $L_\omega$, but not split in any nontrivial extension in $L_\delta/L_\omega$ if $m\geq 1$.  Such places $w$ have density 
\begin{equation}\label{wdensity}
\left\{\begin{aligned}
&\frac{p^m-p^{m-1}}{[L_\delta:F]},\quad & m\geq 1;\\
&\frac{1}{[L_\delta:F]},\quad & m=0.
\end{aligned}\right.
\end{equation}
Fix such a $w$. Since $\omega\left(\CO_{\CB(\delta)}^\times\right)=\mu_{p^m}$, we have
\[\CO_{\CB(\delta)}^\times/U_w=\CO_{\CB(\delta)}^\times/U_\omega\cong \mu_{p^m} \text{ and thus } \CO_w^\times / \CO_w^{\times p^n}\CO_{\CB(\delta)}^\times \cong \BZ/p^{n-m}\BZ.\]
By the identification of characters unramified outside $\CB(\delta)$ and Lemma \ref{character},
\[\ell_{\delta, X}^{-1}(\omega)=\left\{(\omega,\chi_w): \chi_w(\CO_w^\times)=\mu_{p^n} \text{ and }\chi_w(t)=\omega^{-1}(t) \text{ for all }t\in \CO_{\CB(\delta)}^\times\right\},\]
where $\chi_w$ is a character of $\CO_w^\times$, and therefore
\[|\ell_{\delta, X}^{-1}(\omega)\cap \Delta(\delta w)|=\left\{\begin{aligned}
&p^{n-m},\quad & m\geq 1;\\
&p^n-p^{n-1},\quad & m=0.
\end{aligned}\right.\]

For any $\omega\in \Omega_\delta$, 
\[|\ell_{\delta, X}^{-1}(\omega)|=\sum_{\begin{subarray}{c}w\in \CP_1(X)\backslash \ov{\delta}\\U_w=U_\omega\end{subarray}}|\ell_{\delta, X}^{-1}(\omega)\cap \Delta(\delta w)|.\]
Thus by (\ref{wdensity}), Theorem \ref{chebotarev} and Definition \ref{new-convergence}, for $Y>\N_\delta$ and $X>\CL'(Y)$ one has
\[\left|\frac{|\ell_{\delta, X}^{-1}(\omega)|}{|\CP_1(X)\backslash \ov{\delta}|}-\frac{p^{n}-p^{n-1}}{[L_\delta:F]}\right|\leq \frac{p^n}{Y}.\]
\end{proof}

For $\gamma\in \Omega_1$, denote $\Delta^\gamma(\delta*\CP_1(X))=\ell_{\delta, X}^{-1}(\Omega_\delta^\gamma)$.
\begin{prop}\label{global-local2}
For $Y> \max(\N_\delta, 2d(\CB,k))$ and $X>\CL'(Y)$, the map $\ell_{\delta,X}: \Delta^\gamma(\delta*\CP_1(X))\ra \Omega^\gamma_\delta$ is surjective and for any $\omega\in \Omega^\gamma_\delta$, 
\[\left|\frac{|\ell_{\delta, X}^{-1}(\omega)|}{|\Delta^\gamma(\delta*\CP_1(X))|}-\frac{1}{|\Omega^\gamma_\delta|}\right|\leq \frac{8d(\CB,k)}{|\Omega^\gamma_\delta|Y}.\]
\end{prop}
\begin{proof}
Denote
\[A=\frac{p^n-p^{n-1}}{d(\CB,k)},\quad B=|\CP_1(X)\backslash \ov{\delta}|.\]
By Proposition \ref{fiber}, for any $\omega\in \Omega_\delta^\gamma$
\[|\ell_{\delta, X}^{-1}(\omega)|=AB+\theta_\omega B p^n/Y, \text{ with }|\theta_\omega|\leq 1,\]
and
\begin{eqnarray*}
|\Delta^\gamma(\delta*\CP_1(X))|&=&\sum_{\omega'\in \Omega^\gamma_\delta}|\ell_{\delta, X}^{-1}(\omega')|\\
&=&AB|\Omega^\gamma_\delta|+\left(\sum_{\omega'\in \Omega^\gamma_\delta}\theta_{\omega'}\right)Bp^n/Y.
\end{eqnarray*}
Then,
\begin{eqnarray*}
\left|\frac{|\ell_{\delta, X}^{-1}(\omega)|}{|\Delta^\gamma(\delta*\CP_1(X))|}-\frac{1}{|\Omega_\delta^\gamma|}\right|&=&\left|\frac{A+\theta_{\omega}  p^n/Y}{A|\Omega_\delta^\gamma|+\left(\sum_{\omega'\in \Omega_\delta^\gamma}\theta_{\omega'}\right)p^n/Y}-\frac{1}{|\Omega_\delta^\gamma|}\right|\\
&=&\left|\frac{\left(\theta_\omega |\Omega_\delta^\gamma|-\sum_{\omega'\in \Omega_\delta^\gamma}\theta_{\omega'}\right)p^n/Y}{\left(A|\Omega_\delta^\gamma|+\left(\sum_{\omega'\in \Omega_\delta^\gamma}\theta_{\omega'}\right)p^n/Y\right)|\Omega_\delta^\gamma|}\right|\\
&\leq &\left|\frac{2 p^n/Y}{\left(A|\Omega_\delta^\gamma|+\left(\sum_{\omega'\in \Omega_\delta^\gamma}\theta_{\omega'}\right)p^n/Y\right)}\right|\\
&\leq &\left|\frac{2 p^n/Y}{\left(A-p^n/Y\right)|\Omega_\delta^\gamma|}\right|\leq \frac{4p^n}{A|\Omega_\delta^\gamma|Y}.
\end{eqnarray*}
The last inequality follows for $Y\geq 2d(\CB,k)$, whence $p^n/Y\leq A/2$.
\end{proof}

\section{The  distribution of ranks for global Selmer groups}
The bridge built between global and local characters indeed preserves the Selmer ranks for global and local characters (see Proposition \ref{dim-preserving}), and we will transfer the statistical behavior of Selmer groups associated to the local characters to those of the global ones.

\begin{lem}\label{torsion-gamma}
\begin{itemize}
\item[(1)] For  $\chi\in \CC(F)$, $\dim_k \kappa_{A_\chi,\fp}(A_\chi (F)_\tor) =1$ resp. $0$ if $F=F_1$ resp. $F\neq F_1$.
\item[(2)] Suppose $\delta\in \CD$ and $\chi\in \Delta(\delta)$. If $F=F_d$, then $\loc_v(\kappa_{A_\chi,\fp}(A_\chi (F)_\tor)\neq 0$ for all $v\mid \delta$.
\item[(3)] If $F_1=F_d$, then for any $\gamma\in \Omega_1$, $\Sha_{\CP_1}(\gamma)=0$.
\item[(4)] Suppose $\delta\in \CD$ and $\chi\in \Delta(\delta)$ restricting to $\gamma\in \Omega_1$. If $F\neq F_d$ , then $\kappa_{A_\chi,\fp}(A_\chi (F)_\tor)\subset \Sha_{\CP_1}(\gamma)$. Moreover, if $F=F_{d-1}\neq F_d$ with $d\geq 2$ $($necessarily $(\CMS)$ holds$)$, then  $\Sha_{\CP_1}(\gamma)=\kappa_{A_\chi,\fp}(A_\chi (F)_\tor)$ has dimension one.
\end{itemize}
\end{lem}
\begin{proof}
\noindent (1) Through the Kummer $\fp$-descent map
\[\kappa_{A_\chi,\fp}(A_\chi(F)_\tor)=\kappa_{A_\chi,\fp}(A_\chi[\fp^\infty](F))\cong A_\chi[\fp^\infty](F)\otimes_\CO k.\]
By Theorem \ref{twist-av}, $A_\chi[\fp](F)=A[\bp](F)$, and therefore (1) follows.  

\noindent (2) Suppose $\delta\in \CD$ and $\chi\in \Delta(\delta)$. By Theorem \ref{twist-av}, $A_\chi[\fp^d]\cong A[\bp^d]$ as $G_F$-modules. If $F=F_d$, then $A_\chi[\fp^{d}](F)=A_\chi[\fp^d]$. For $v\mid \delta$, $\chi_v$ is totally ramified of order $p^n$. By Lemma \ref{gen}, $A_\chi[\fp^{d+1}](F_v)=A_\chi[\fp^d]$, and thus $A_\chi[\fp^d]$ has nonzero image  in $\RH^1(F_v,A_\chi[\fp])$ under the Kummer map.

\noindent (3) Suppose $F_1=F_d$. For any $\gamma\in \Omega_1$, $\delta\in \CD$ and $\omega\in \Omega_\delta^\gamma$, $F_{\Sha_{\CP_1}(\gamma)}=F_\omega\cap F_d=F_1$ and thus $\Sha_{\CP_1}(\gamma)=0$ by Proposition \ref{disjoint}.

\noindent (4) Suppose $F\neq F_d$. Let $0\leq i\leq d-1$ be the integer such that $F=F_i$ but $F\neq F_{i+1}$. We have $A_\chi [\fp^d]\cong A[\bp^d]$ as $G_F$-modules. If $i=0$, then $\kappa_{A_\chi, \fp}(A_\chi[\fp^{i}])=0$. Assume $i>0$. Then  $\kappa_{A_\chi, \fp}(A_\chi[\fp^{i}])\neq 0$. For any $v\in \CP_1$, $A_\chi [\fp^d]\subset A_\chi(F_v)$, and thus $\loc_v(\kappa_{A_\chi,\fp}(A_\chi[\fp^i]))=0$. Therefore $\kappa_{A_\chi,\fp}(A_\chi[\fp^i])\subset \Sha_{\CP_1}(\gamma)$. Moreover suppose $F=F_{d-1}\neq F_d$ with $d\geq 2$. Necessarily $(\CMS)$ holds and then through the character $\theta$ of (\ref{CM-character}),
\[\Gal(F_d/F_{d-1})\cong (1+\bp^{d-1})/(1+\bp^d)\cong \BF_p,\]
i.e. $F_d/F_1$ is cyclic of order $p$. Since $F=F_{d-1}$, $\kappa_{A_\chi,\fp}(A_\chi[\fp^i])\neq 0$. Thus we must have $F_{\Sha_{\CP_1}(\gamma)}=F_d$ and $\Sha_{\CP_1}(\gamma)=\kappa_{A_\chi,\fp}(A_\chi (F)_\tor)$ of dimension one by Proposition \ref{disjoint}.
\end{proof}


For $\chi\in \Delta(F)$, denote $r(\chi)=\dim_{k} \Sel_{\fp}(A_\chi).$

\begin{prop}\label{dim-preserving}
\begin{itemize}
\item[1.] Suppose  $F\neq F_d$,   $\delta\in \CD$ and $\ell_{\delta,X}: \Delta(\delta,X)\ra \Omega_\delta$ is the restriction map. For any $\chi\in \Delta(\delta,X)$, $r(\chi)=r(\ell_{\delta,X}(\chi))$.
\item[2.] Suppose $F=F_d$,  $\delta\in \CD$ and $\ell_{\delta,X}: \Delta(\delta*\CP_1(X))\ra \Omega_\delta$ is the restriction map. For any $\chi\in \Delta(\delta*\CP_1(X))$, $r(\chi)=r(\ell_{\delta,X}(\chi))+1$.
\end{itemize}
\end{prop}
\begin{proof}
The first assertion follows immediately from Proposition \ref{Selmer}. Suppose $F=F_d$, $\delta\in \CD$, $w\in \CP_1\backslash \ov{\delta}$ and $\chi\in \Delta(\delta w)$. For $\omega'=(\chi_v)_{v\in \CB(\delta w)}$, by Proposition \ref{Selmer},
\[\Sel_{\fp}(A_\chi)=\Sel(T,\omega').\]
Let $\omega=(\chi_v)_{v\in \CB(\delta)}=\ell_{\delta, X}(\chi)$, and then $\omega'\in \eta_{\delta,w}^{-1}(\omega)$. By Lemma \ref{torsion-gamma}.(2), we have $\loc_w(\Sel(T,\omega'))=\loc_w(\Sel_{\fp}(A_\chi))\neq 0$. Hence by the proof of Proposition \ref{selmer-change}, $r(\omega')-r(\omega)=1$ and the second assertion follows as desired.
\end{proof}

\begin{defn}\label{stratification0}
For $\gamma\in \Omega_1$, $\delta\in \CD$ and $X>\CL(\N_\delta)$ resp. $\CL'(\N_{\delta})$,  let $\ell_{\delta,X}$ be the map in Proposition \ref{global-local} resp. \ref{fiber}, and define
\begin{equation*}
\ell^{-1}_{\delta,X}(\Omega^\gamma_\delta)=:\left\{\begin{aligned}
&\Delta^\gamma(\delta, X),\quad&&\text{if $F\neq F_d$};\\
&\Delta^\gamma(\delta*\CP_1(X)),\quad&&\text{if $F=F_d$}.
\end{aligned}\right.
\end{equation*}
\end{defn}
If $F\neq F_d$, $\chi \in \Delta^\gamma(\delta,X)$, then $r(\chi)\geq r_\gamma$; If  $F= F_d$, $\chi \in \Delta^\gamma(\delta*\CP_1(X))$, then $r(\chi)\geq r_\gamma+ 1=1$. In the latter case, by Lemma \ref{torsion-gamma}.(3), $r_\gamma=0$.

\begin{prop}\label{passage}
\begin{itemize}
\item[(1)] Suppose $F\neq F_d$, $\gamma\in \Omega_1$, $\delta\in \CD$ and $r\geq r_\gamma$. If $X> \CL(\N_\delta)$, then
\[\frac{|\{\chi\in \Delta^\gamma(\delta,X):  r(\chi)=r\}|}{|\Delta^\gamma(\delta,X)|}=E^\gamma_\delta(r-r_\gamma).\]
\item[(2)] Suppose $F=F_d$, $\gamma\in \Omega_1$, $\delta\in \CD$ and $r\geq 1$. If  $Y>\max(\N_{\delta}, 2d(\CB,k))$ and $X>\CL'(Y)$, then
\[\left|\frac{|\{\chi\in \Delta^\gamma(\delta*\CP_1(X)):  r(\chi)=r\}|}{|\Delta^\gamma(\delta*\CP_1(X))|}-E^\gamma_\delta(r-1)\right|\leq \frac{8d(\CB,k)}{Y}.\]
\end{itemize}
\end{prop}
\begin{proof}
First suppose $F\neq F_d$, $\gamma\in \Omega_1$, $\delta\in \CD$ and $r\geq r_\gamma$. By Proposition \ref{global-local}, the map $\ell_{\delta,X}:\Delta^\gamma(\delta,X)\ra\Omega^\gamma_\delta$ is surjective and all fibers have the same cardinality. By Proposition \ref{dim-preserving}, $r(\chi)=r(\ell_{\delta,X}(\chi))$ and then
\[\frac{|\{\chi\in \Delta^\gamma(\delta,X):  r(\chi)=r\}|}{|\Delta^\gamma(\delta,X)|}=\frac{|\{\omega\in \Omega^\gamma_\delta:  r(\omega)=r\}|}{|\Omega^\gamma_\delta|}=E^\gamma_\delta(r-r_\gamma).\]

Next suppose $F=F_d$, $\gamma\in \Omega_1$, $\delta\in \CD$ and $r\geq 1$.  By Proposition \ref{dim-preserving}, $r(\chi)=r(\ell_{\delta,X}(\chi))+1$, and then
\begin{equation}\label{eq}8d(\CB,k)
\frac{|\{\chi\in \Delta^\gamma(\delta*\CP_1(X)):r(\chi)=r\}|}{|\Delta^\gamma(\delta*\CP_1(X))|}
=\sum_{\begin{subarray}{c}\omega\in \Omega^\gamma_\delta\\r(\omega)=r-1\end{subarray}}\frac{|\ell_{\delta, X}^{-1}(\omega)|}{|\Delta^\gamma(\delta*\CP_1(X))|}.
\end{equation}
By Proposition \ref{global-local2}, for $Y\geq \max(\N_{\delta},2d(\CB,k))$ and $X>\CL'(Y)$ and every $\omega\in \Omega_\delta^\gamma$
\[\left|\frac{|\ell_{\delta, X}^{-1}(\omega)|}{|\Delta^\gamma(\delta*\CP_1(X))|}-\frac{1}{|\Omega^\gamma_\delta|}\right|\leq \frac{8d(\CB,k)}{|\Omega^\gamma_\delta|Y}.\]
By (\ref{eq}) and the above inequality, 
\begin{eqnarray*}
&&\left|\frac{|\{\chi\in \Delta^\gamma(\delta*\CP_1(X)):  r(\chi)=r\}|}{|\Delta^\gamma(\delta*\CP_1(X))|}-E^\gamma_\delta(r-1)\right|\\
&=&\left|\sum_{\begin{subarray}{c}\omega\in \Omega^\gamma_\delta\\r(\omega)=r-1\end{subarray}}\left(\frac{|\ell_{\delta, X}^{-1}(\omega)|}{|\Delta^\gamma(\delta*\CP_1(X))|}-\frac{1}{|\Omega^\gamma_\delta|}\right) \right| \\
&\leq &\sum_{\begin{subarray}{c}\omega\in \Omega^\gamma_\delta\\r(\omega)=r-1\end{subarray}}\left|\frac{|\ell_{\delta, X}^{-1}(\omega)|}{|\Delta^\gamma(\delta*\CP_1(X))|}-\frac{1}{|\Omega^\gamma_\delta|}\right| \leq \frac{8d(\CB,k)}{Y}.
\end{eqnarray*}
\end{proof}

\begin{defn}\label{stratification}
\begin{itemize}
\item[(1)] For $\gamma\in \Omega_1$ and integers $k\geq 0$, define
\begin{equation*}
\Delta^\gamma_k(X):=\left\{\begin{aligned}
&\bigsqcup_{\delta\in \CD_k(X)}\Delta^\gamma(\delta, \CL(L_{k+1}(X))),\quad&&\text{if  $F\neq F_d$};\\
&\bigsqcup_{\delta\in \CD_{k}(X)}\Delta^\gamma(\delta*\CP_1(\CL'(L_{k+1}(X)))),\quad&&\text{if $F=F_d$}.
\end{aligned}\right.
\end{equation*}
Moreover, for $k=-1$, define $\Delta^\gamma_{-1}(X)=\emptyset$ resp. $\Delta^{\gamma}(1)$ if $F\neq F_d$ resp. $F=F_d$, where $\Delta^\gamma(1)$ denotes the set of characters ramified only in $\CB$ and restricting to $\gamma$ over $\CB$.
\item[(2)] The union $\Delta^\gamma(F)=\bigcup_{k,X} \Delta_k^\gamma(X)$ is the set of characters $\chi\in \Delta(F)$ which restrict to $\gamma$ over $\CB$ and we call the collection of $\Delta_k(X)$ a ``fan-structure"  stratification  on $\Delta^\gamma(F)$.
\end{itemize}
\end{defn}

\begin{lem}\label{samecard}
Suppose $\gamma\in \Omega_1$ and $F\neq F_d$. For $\delta\in \CD_k$  and $X>\CL(\N_\delta)$, $|\Delta^\gamma (\delta,X)|=\varphi(p^n)^k|\Delta^\gamma(1,X)|$. In particular 
\[|\Delta_k^\gamma(X)|=|\Delta^\gamma(\delta,\CL(L_{k+1}(X)))||\CD_k(X)|.\]
\end{lem}
\begin{proof}
Suppose $\gamma\in \Omega_1$ and $\delta=v_1v_2\cdots v_k$. For $1\leq i\leq k$, denote $\omega_i= (\gamma'=1,\mu_i)\in \Omega^1_{v_i}$ where $\mu_i$ is a fixed character in $\CC_n^\tr(F_{v_i})$. For $X>\CL(\N_\delta)$, by Proposition \ref{global-local}, there exists $\chi_i\in \Delta^1(v_i, X)$  such that $\ell_{v_i,X}(\chi_i)=\omega_i$ for each $i$. Then there exist a bijection
\[((\BZ/p^n\BZ)^\times)^k\times \Delta^\gamma(1,X)\ra \Delta^\gamma(\delta,X)\]
given by $(n_1,\cdots,n_k,\psi)\mapsto \chi_1^{n_1}\cdots \chi_k^{n_k}\psi$. The lemma follows as desired.
\end{proof}

\begin{lem}\label{samecard1}
Suppose $\gamma\in \Omega_1$ and $F=F_d$. For $\delta\in \CD_k(X)$,
\[\lim_{X\ra \infty}\frac{|\CD_k(X)||\Delta^\gamma(\delta*\CP_1(\CL'(L_{k+1}(X))))|}{|\Delta_k^\gamma(X)|}=1.\]
\end{lem}
\begin{proof}
Denote $Y=L_{k+1}(X)$, $X'=\CL'(Y)$ and $A=(p^n-p^{n-1})/d(\CB,k)$.  For $\delta\in \CD_k(X)$, put $B=|\Omega_\delta^\gamma|$.  For $\omega\in \Omega^\gamma_\delta$, by Proposition \ref{fiber} and \ref{global-local2}, we have 
\[|\ell^{-1}_{\delta,X'}(\omega)|=|\CP_1(X')\backslash \ov{\delta}|A(1+\alpha_{\delta,\omega}p^n/AY),\]
\[B|\ell^{-1}_{\delta,X'}(\omega)|=|\Delta^\gamma(\delta*\CP_1(X'))|(1+\beta_{\delta,\omega}8d(\CB,k)/Y)\]
where $\alpha_{\delta,\omega}$ and $\beta_{\delta,\omega}$ depend on $X$ and have absolute values $\leq 1$. For $X>16d(\CB,k)$ so that $|8\beta_{\delta,\omega}d(\CB,k)/Y|<1/2$, there exists an absolute constant $C_0>0$ and $c_{\delta}$ for every $\delta$ with $|c_{\delta}|\leq C_0$ such that 
\begin{equation}\label{fiber-card}
|\Delta^\gamma(\delta*\CP_1(X'))|=AB|\CP_1(X')\backslash \ov{\delta}|(1+c_{\delta}/Y).
\end{equation}

Denote $\omega(\ov{\delta})$ to be the number of primes dividing $\ov{\delta}$ and let $\pi(X)$ be the number of primes of $F$ of norm $\leq X$. By the generalized prime number theorem (cf. \cite{LO77}), as $X\ra\infty$, 
\begin{equation}\label{pnt}
\pi(X)=X/\log X+O(X/(\log X)^2).
\end{equation}
Let $X$ be sufficiently large so that $|\CP_1(X')|=\pi(X')-|\CB|$ and $C_0/Y<1/2$. By definition, $L_{k+1}(X)\geq XL_k(X)$. Thus for $\delta'\in \CD_k(X)$, using the estimation (\ref{pnt}), as $X\ra \infty$,
\[\omega(\ov{\delta'})\leq \pi(L_{k}(X))=o(\pi(L_{k+1}(X)))=o(\CP_1(X')).\]
Then, by (\ref{fiber-card}),
\begin{eqnarray*}
&&\frac{|\CD_k(X)||\Delta^\gamma(\delta*\CP_1(\CL'(L_{k+1}(X))))|}{|\Delta_k^\gamma(X)|}\\
&=&\frac{|\CD_k(X)||\CP_1(X')|+|\CD_k(X)||\CP_1(X')|c_\delta/Y-|\CD_k(X)|\omega(\ov{\delta})(1+c_{\delta}/Y)}{|\CD_k(X)||\CP_1(X')|+ |\CP_1(X')|\sum_{\delta'} c_{\delta'}/Y-\sum_{\delta'}\omega(\ov{\delta'})(1+c_{\delta'}/Y)}
\end{eqnarray*}
converges to $1$ as $X\ra \infty$.


\end{proof}

\begin{thm}\label{asymptotic-markov}
Suppose $\gamma\in \Omega_1$. For $r\geq r_\gamma \text{ resp. } 1$ if $F\neq F_d$ resp. $F=F_d$, one has
\[\lim_{X\ra \infty} \frac{|\{\chi\in \Delta^\gamma_k(X):r(\chi)=r\}|}{|\Delta^\gamma_k(X)|}=\left\{\begin{aligned}
&(E_1^\gamma M_\epsilon^k)(r-r_\gamma),\quad&&\text{if  $F\neq F_d$};\\
&(E_1^\gamma M_\epsilon^k)(r-1),\quad&&\text{if $F=F_d$}.
\end{aligned}\right..\]
\end{thm}
\begin{proof}
If $F\neq F_d$ resp. $F=F_d$, denote $r_0=r_\gamma$ resp. $1$ and $\Delta_{\gamma,\delta,X}=\Delta^\gamma(\delta,\CL(L_{k+1}(X)))$ resp. $\Delta^\gamma(\delta*\CP_1(\CL'(L_{k+1}(X))))$.  By Proposition \ref{passage}, for sufficiently large $X$ and $Y=L_{k+1}(X)$,
\begin{eqnarray*}
&&\left|\frac{|\{\chi\in \Delta^\gamma_k(X):r(\chi)=r\}|}{|\Delta^\gamma_k(X)|}-E^\gamma_{\CD_k(X)}(r-r_0)\right|\\
&=&\left|\sum_{\delta\in \CD_k(X)} \frac{|\chi\in \Delta_{\gamma,\delta,X}:r(\chi)=r|}{|\Delta_{\gamma,\delta,X}|}\cdot \frac{|\Delta_{\gamma,\delta,X}|}{|\Delta^\gamma_k(X)|}-\frac{E^\gamma_\delta(r-r_0)}{|\CD_k(X)|}\right|\\
&=&\left|\sum_{\delta\in \CD_k(X)} \left(\frac{|\chi\in \Delta_{\gamma,\delta,X}:r(\chi)=r|}{|\Delta_{\gamma,\delta,X}|}-E_\delta^\gamma(r-r_0)\right) \frac{|\Delta_{\gamma,\delta,X}|}{|\Delta^\gamma_k(X)|}+E^\gamma_\delta(r-r_0)\left(\frac{|\Delta_{\gamma,\delta,X}|}{|\Delta^\gamma_k(X)|}-\frac{1}{|\CD_k(X)|}\right)\right|\\
&\leq& \frac{8d(\CB,k)}{Y}\frac{|\CD_k(X)||\Delta_{\gamma,\delta,X}|}{|\Delta^\gamma_k(X)|}+\left|\frac{|\Delta_{\gamma,\delta,X}|}{|\Delta^\gamma_k(X)|}-\frac{1}{|\CD_k(X)|}\right|
\end{eqnarray*}
which converges to zero as $X\ra \infty$ by Lemma \ref{samecard} and \ref{samecard1}. By Theorem \ref{local-markov}, $E^\gamma_{\CD_k(X)}(r-r_0)$ converges to $(E^\gamma_1M_\epsilon^k)(r-r_0)$ as $X\ra\infty$.
\end{proof}

\begin{thm}\label{main-body}
Suppose $\gamma\in \Omega_1$. Under the ``fan-structure" stratification $\Delta^\gamma(F)=\bigcup_{k,X} \Delta_k^\gamma(X)$,  for $r\geq r_\gamma \text{ resp. } 1$,  one has
\[\lim_{k\ra \infty}\lim_{X\ra \infty} \frac{|\{\chi\in \Delta^\gamma_k(X):r(\chi)=r\}|}{|\Delta^\gamma_k(X)|}=\left\{\begin{aligned}
&\sD_q^\epsilon(r-r_\gamma),\quad&&\text{if  $F\neq F_d$},\\
&\sD_q^\epsilon(r-1),\quad&&\text{if $F=F_d$}.
\end{aligned}\right.\]
Here $\epsilon=1$ or $1/2$ according to  whether $(\CMS)$ or $(\CMU)$ holds.
\end{thm}
\begin{proof}
This is an immediate consequence of Theorem \ref{local-markov} and \ref{asymptotic-markov}.
\end{proof}
 If we denote $\Delta_k(X)=\bigcup_{\gamma\in \Omega_1}\Delta^\gamma_k(X)$, then $\Delta(F)=\bigcup_{k,X} \Delta_k(X)$.
\begin{coro}\label{main2-body}
Assume either $F=F_{d-1}\neq F_d$ with $d\geq 2$ or $F_1=F_d$. Then
\[\lim_{k\ra \infty}\lim_{X\ra \infty} \frac{|\{\chi\in \Delta_k(X):\dim_k \left(\Sel_\fp(A_\chi)/\kappa_{A_\chi,\fp}(A_\chi(F)_\tor)\right)=r\}|}{|\Delta_k(X)|}=\sD_q^\epsilon(r).\]
\end{coro}
\begin{proof}
It follows from Lemma \ref{torsion-gamma} and Theorem \ref{main-body} that the distribution of $\Sel_\fp(A_\chi)$ with $\chi\in \Delta^\gamma(F)$ is described in term of $\sD_q^\epsilon$, shifted by the dimension of $(\kappa_{A_\chi,\fp}(A_\chi(F)_\tor)+\Sha_{\CP_1}(\gamma))$. Under the hypothesis of the corollary, for any $\gamma\in \Omega_1$, $\Sha_{\CP_1}(\gamma)\subset \kappa_{A_\chi,\fp}(A_\chi(F)_\tor)$. Thus the corollary follows.
\end{proof}

The ``fan-structrue" stratification $\Delta^\gamma(F)=\bigcup_{k,X} \Delta_k^\gamma(X)$ is in fact a disjoint union for $k$ and $\Delta_k^\gamma(X)$ consists of characters ramified exactly at $k$ (resp. $k+1$) places of $\CP_1$ if $F\neq F_d$ (resp. $F=F_d$). If we denote $\Delta^\gamma(k,X)=\bigcup_{i= -1}^k \Delta_i^\gamma(X)$, then
\[\lim_{k\ra \infty}\lim_{X\ra \infty}\Delta^\gamma(k,X)=\Delta^\gamma(F).\]
Therefore one may prefer a limit along $\Delta^\gamma(k,X)$. Since characters ramified at more places have much larger cardinality, $\Delta_k^\gamma(X)$ dominants the union $\Delta^\gamma(k,X)$ and the limit doesn't change. Therefore we have the following variants of Theorem \ref{main-body} and Corollary \ref{main2-body}.
\begin{thm}\label{main-body1}
Suppose $\gamma\in \Omega_1$. Under the ``fan-structure" stratification,  for $r\geq r_\gamma \text{ resp. } 1$,  
\[\lim_{k\ra \infty}\lim_{X\ra \infty} \frac{\left|\{\chi\in \Delta^\gamma(k,X):r(\chi)=r\}\right|}{\left|\Delta^\gamma(k,X)\right|}=\left\{\begin{aligned}
&\sD_q^\epsilon(r-r_\gamma),\quad&&\text{if  $F\neq F_d$},\\
&\sD_q^\epsilon(r-1),\quad&&\text{if $F=F_d$}.
\end{aligned}\right.\]
\end{thm}
Denote $\Delta(k,X)=\bigcup_{\gamma\in \Omega_1} \Delta^\gamma(k,X)$.
\begin{coro}\label{main2-body1}
Assume either $F=F_{d-1}\neq F_d$ with $d\geq 2$ or $F_1=F_d$. Then
\[\lim_{k\ra \infty}\lim_{X\ra \infty} \frac{\left|\{\chi\in \Delta(k,X):\dim_k \left(\Sel_\fp(A_\chi)/\kappa_{A_\chi,\fp}(A_\chi(F)_\tor)\right)=r\}\right|}{\left| \Delta(k,X)\right|}=\sD_q^\epsilon(r).\]
\end{coro}
\begin{proof}[Proof of Theorem \ref{main-body1}]
For $k\geq 0$, for simplicity, denote $\Delta_{\gamma, \delta,X}=\Delta^\gamma(\delta, \CL(L_{k+1}(X)))$ resp. $\Delta^\gamma(\delta*\CP_1(\CL'(L_{k+1}(X))))$. Then
\[\Delta^\gamma_k(X)=\bigsqcup_{\delta\in \CD_k(X)} \Delta_{\gamma,\delta,X}.\]
For $k\geq -1$, let $S_k(X)$ be the number of $\chi\in \Delta^\gamma_k(X)$ such that $r(\chi)=r$ and denote $\alpha_k(X)=S_k(X)/|\Delta_k^\gamma(X)|$. Note $\alpha_{-1}(X)$ is a constant. For $k\geq 0$, by Theorem \ref{asymptotic-markov} and \ref{main-body},
\[\lim_{X\ra\infty} \alpha_k(X)=\alpha_k \text{ and }  \lim_{k\ra\infty} \alpha_k=\alpha,\]
where $\alpha_k$ and $\alpha$ denote the limits as indicated in the theorems. Thus for any $\epsilon>0$, there exists $k_0$ such that for $i> k_0$, $|\alpha_i-\alpha|<\epsilon$. Let $k>k_0$ be a large integer. There exist $X_0$ and $C_\epsilon$ such that for $X>X_0$, 
\begin{eqnarray*}
\left|\frac{\sum_{i=-1}^{k}S_k(X)}{\sum_{i=-1}^k |\Delta^\gamma_k(X)|}-\alpha\right|&=&\left|\frac{\sum_{i=-1}^{k_0} (\alpha_i(X)-\alpha)|\Delta^\gamma_i(X)|+\sum_{i=k_0+1}^{k} (\alpha_i(X)-\alpha)|\Delta^\gamma_i(X)|}{\sum_{i=-1}^k|\Delta^\gamma_i(X)|}\right|\\
&\leq &\frac{\sum_{i=-1}^{k_0} C_\epsilon |\Delta^\gamma_i(X)|+2\epsilon \sum_{i=k_0+1}^{k} |\Delta^\gamma_i(X)|}{\sum_{i=-1}^k|\Delta^\gamma_i(X)|}.
\end{eqnarray*}
Now it suffices to prove for this fixed $k$, as $X\ra \infty$, 
\begin{equation}\label{0-density}
\frac{\sum_{i=-1}^{k_0}  |\Delta^\gamma_i(X)|}{\sum_{i=-1}^k|\Delta^\gamma_i(X)|}\ra 0.
\end{equation}
For any $0\leq i\leq k_0$, by Lemma \ref{samecard} and \ref{samecard1} and (\ref{fiber-card}), for large $X$,
\[\frac{|\Delta^\gamma_i(X)|}{|\Delta^\gamma_k(X)|}\leq \frac{|\CD_i(X)|}{|\CD_k(X)|}\]
which converges to zero by Lemma \ref{0-density1}. If $i=-1$, then $|\Delta^\gamma_{-1}(X)|$ is a constant and $|\Delta^\gamma_{-1}(X)|/|\Delta^\gamma_k(X)|$ also converges to $0$. Thus (\ref{0-density}), and hence the theorem, follow as desired. 
\end{proof}

\begin{lem}\label{0-density1}
For $k\geq 0$,
\[\lim_{X \ra \infty}\frac{|\CD_k(X)|}{|\CD_{k+1}(X)|}=0.\]
\end{lem}
\begin{proof}
Let $X$ be sufficiently large so that any $v\in \CB$ has norm $<X$. By definition $L_{k+1}(X)\geq XL_k(X)$ for $k\geq 1$ and we accept $L_0(X)=1$. For any $\delta=v_1\cdots v_k\in \CD_k(X)$ and any $v_{k+1}\notin \CB$ with norm in  $[L_k(X),L_{k+1}(X))$, $\delta v_{k+1}\in \CD_{k+1}(X)$. If we denote $n_k(X)$ to be the number of places not in $\CB$ and of norm in $[L_k(X),L_{k+1}(X))$, then $|\CD_{k+1}(X)|\geq |\CD_k(X)|n_k(X)$. Let $\pi(X)$ be the number of places of norm $<X$.  If $k=0$, then $n_0(X)=\pi(L_1(X))-|\CB|$. For $k\geq 1$,
\begin{eqnarray*}
n_k(X)=\pi(L_{k+1}(X))-\pi(L_k(X))\geq \pi(XL_k(X))-\pi(L_k(X)).
\end{eqnarray*}
Thus by the estimation (\ref{pnt}), as $X\ra \infty$, $n_k(X)\ra \infty$ as desired. 
\end{proof}

\section{Supplementary proofs}
\begin{proof}[Proof of Theorem \ref{growth}]
For each $L\in \Xi$, let $m_L$ be the number of characters in $\pi^{-1}(L)$ such that $\Sel_{\fp}(A_\chi)/\kappa_{A_\chi,\fp}(A_\chi(F)_\tor)=0$. Through the decomposition (\ref{eigen-decomposition}), we see $\rank_\BZ A(L)=\rank_\BZ A(F)$ if $m_L=p-1$.
Let $A_{k}(X)$ be the number of $L\in\Xi(k,X)= \bigcup_{i\geq -1}^k\Xi_i(X)$  such that $m_L=p-1$. By Theorem \ref{main}, as $X\ra\infty$ and then $k\ra \infty$,
\[\frac{A_k(X)}{|\Xi(k,X)|}+\left(\frac{p-2}{p-1}\right)\left(1-\frac{A_{k}(X)}{|\Xi(k,X)|}\right)\geq    \frac{\sum_{L\in \Xi(k,X)}m_L}{(p-1)|\Xi(k,X)|}\ra \sD_q^\epsilon(0).\]
Then by Proposition \ref{distribution}, 
\begin{eqnarray*}
\frac{A_k(X)}{|\Xi(k,X)|}&\geq & (p-1)\left(\frac{\sum_{L\in \Xi(k,X)}m_L}{(p-1)|\Xi(k,X)|}-\frac{p-2}{p-1} \right)\\
&\ra &(p-1)\left(\sD_q^\epsilon(0)-\frac{p-2}{p-1}\right)\\
&>& (p-1)\left(\frac{1}{p-1}-\frac{q^\epsilon}{q-1}\right)\geq0.
\end{eqnarray*}
If $p=2$, we obtain $\sD_q^\epsilon(0)$ as a lower bound. This proves Theorem \ref{growth}. 
\end{proof}

\begin{proof}[Proof of Theorem \ref{growth1}]
For $L\in \Xi^1$, we first show that $\Sel_\fp(A_\chi)$ are all isomorphic for all $\chi\in \pi^{-1}(L)$. This is done by comparing local Kummer images everywhere. Let $\chi$ be a character of $\pi^{-1}(L)$. For $v\in \CB$,  $\chi$ restricts to the trivial character at $v$ and $A_\chi=A_1$ over $F_v$. Thus the local Kummer image $\CL_{A,v,\chi_v}$ coincides with that of $A_1/F_v$. For $v\in \CP_0$, or $v\in \CP_1$ and $\chi_v$ is unramified, by Proposition \ref{local-condition1}, $\CL_{A,v,\chi_v}$ is the unramified cohomology. For  $v\in \CP_1$ and $\chi_v$ is ramified, by Proposition \ref{local-condition3}, $\CL_{A,v,\chi_v}$ is determined by the local field $L_v$. Thus $\Sel_\fp(A_\chi)$ only depends on $L$  for $\chi\in \pi^{-1}(L)$. 

Also all the $G_F$-modules $A_\chi[\fp]$ are isomorphic. Thus the number $m_L$ of characters in $\pi^{-1}(L)$ such that $\Sel_{\fp}(A_\chi)/\kappa_{A_\chi,\fp}(A_\chi(F)_\tor)=0$ is either $p-1$ or $0$.  As the argument of the proof of Theorem \ref{growth}, we obtain $\sD_q^\epsilon$(0) as a lower bound for the proportion of $L\in \Xi^1$ without rank growth. 
\end{proof}
\begin{proof}[Proof of Theorem \ref{growth2}]
Under the hypotheses of the theorem, by Appendix \ref{stc}, we may replace the data $(A,\lambda,\BK,\BO,p,\bp)$ by a new one $(A_1,\lambda_1,K,\CO,p,\fp)$ with $\omega\in \CO=\CO_K$ and $p\nmid \deg(\lambda_1)$, if necessary. Then all the twists $A_\chi$ admit CM by $\CO$ and $\dag$-sesquilinear  symmetric isogenies of degree prime to $p$ (see the paragraphs preceding Proposition \ref{id-WP}). We are supposed to be in the symplectic case. Also by assumption, $\Sha(A_\chi)[\fp^\infty]$ is finite. Then by \cite[Corollary 6.4-6.5]{Shu-quadratic}, $\Sha(A_\chi)[\fp]$ has even dimension and 
\[\dim_{k} A_\chi(F)\otimes_{\CO} k\equiv \dim_{k} \Sel_\fp(A_\chi)\mod 2.\]
Thus as seen in the proof of Corollary \ref{main2}, 
\[\rank_\CO A_\chi(F)=\dim_{k} A_\chi(F)/(\fp A_\chi(F)+A_\chi(F)_\tor) \equiv \dim_k \Sel_\fp(A_\chi)/\kappa_{A,\chi}(A_\chi(F)_\tor)\mod 2.\]
By Theorem \ref{main} and Proposition \ref{distribution}, the proportion of $\chi$ for which $A_\chi$ has odd $\CO$-rank  is $(1-\beta^1)/2$. 
This proves (1). 

Thus the proportion of $\chi$ for which $\rank_\BZ A_\chi(F)>0$ is at least $(1-\beta^1)/2$. By assuming the worst situation that for $L\in \Xi$ if  $\Sel_\fp(A_\chi)/\kappa_{A,\chi}(A_\chi(F)_\tor)$ has odd dimension for some $\chi\in \pi^{-1}(L)$, then it is so for all characters in $\pi^{-1}(L)$, we get the lower bound of (2).
\end{proof}

\appendix
\section{}\label{STC}
\subsection{Cyclotomic extensions}\label{cyc-extensions}
Let $p$ be a prime and $\omega$ a primitive $p^n$-th root of unity.  Let $\BK$ be a number field with $\BO\subset \BK$ an order. Let $\CO=\BO[\omega]$ and $K=\BK(\omega)$. Let $\bp$ be an invertible prime ideal of $\BO$ above $p$ and let $\fp=(\bp,1-\omega)$ be the unique prime ideal of $\CO$ above $\bp$. Suppose $\fp$ is an invertible prime ideal of $\CO$ and denote $d=\ord_\fp(1-\omega)$.
\begin{lem}\label{cyc-ext}
\begin{itemize}
\item[(1)] For $1\leq i\leq d$, $\BO/\bp^i\cong\CO/\fp^i$.
\item[(2)] If $\fp$ is ramified above $\bp$, then $d=1$.
\item[(3)] If $\omega\notin \BO_\bp$, then $\fp$ is ramified above $\bp$.
\end{itemize}
\end{lem}
\begin{proof}
(1) By the definition of $d$, for $1\leq i\leq d$, $\fp^i=(\bp^i,1-\omega)$ and thus
\[\CO/\fp^i=\BO[\omega]/(\bp^i,1-\omega)\cong \BO/\bp^i.\]

\noindent (2) Suppose $\fp$ is ramified over $\bp$. In particular $\fp^2\cap \BO=\bp$ and we have an injection $\BO/\bp=\BO/(\fp^2\cap\BO)\hookrightarrow \CO/\fp^2$. It follows from (1) that $d=1$.

\noindent (3) Since both $\bp$ and $\fp$ are invertible, $\BO_\bp[\omega]=\CO_\fp=\CO_{K,\fp}.$
Let $f$ be the monic minimal polynomial of $\omega$ over $\BO_\bp$. By \cite[Chap. III, \S 6, Corollary 2]{Serre-localfields}, 
\[\partial_{\CO_{K,\fp}/\BO_\bp}=(f'(\omega))=\fp^{d(\deg(f)-1)}\CO_{K,\fp}.\]
Thus if $\deg(f)\geq 2$, $\fp$ is ramified over $\bp$. Otherwise $\omega\in \BO_\bp$.
\end{proof}

Suppose $\BK$ is a CM field, $\BO=\CO_\BK$ its ring of integers and $\BK_0$ its maximal real subfield. Then $K=\BK(\omega)$ is also a CM field and denote the complex conjugation on $K$ by $\dag$.  Denote $M=\BK\cap \BQ(\omega)$.
\begin{lem}\label{cyc-ext2}
Assume $p$ is prime to the relative discriminant $d_{\BK/M}$.
\begin{itemize}
\item[(1)] $\CO=\BO[\omega]$ is the ring of integers of $K$.
\item[(2)]  If $\bp$ is inert resp. ramified in $\BK/\BK_0$, then so is $\fp$ in $K/K_0$.
\item[(3)] The different $\partial_{\CO/\BO}$ is principal with a generator $\alpha$ such that 
\[\alpha^{1-\dag}=\left\{\begin{aligned}
&1,&\quad \text{if $p=2$};\\
&(-1)^{[\BQ(\omega):M]-1},&\quad \text{if $p\neq 2$}.
\end{aligned}\right.\]
\end{itemize}
\end{lem}
\begin{proof}
(1) Since $\BQ(\omega)$ is Galois over $M$, $\BQ(\omega)$ and $K$ are linearly disjoint over $M$. Since $p$ is prime to $d_{K/M}$, $d_{\BQ(\omega)/M}$ and $d_{K/M}$ are relatively prime. It follows from \cite[I, Proposition 2.11 and its Remark]{Neu92} that $\CO$ is the ring of integers of $K$.  

\noindent (2) First note $K=K_0\BK$ and $\fp$ is stable under $\dag$. Thus $\fp$ is either inert or ramified in $K/K_0$. Put $\CO_0=K_0\cap \CO$ and $\fp_0=K_0\cap \fp$. If $\bp$ is inert in $\BK/\BK_0$, then $\fp_0$ is unramfied, hence inert, in $K/K_0$. If $\bp$ is ramified in $\BK/\BK_0$, then the residue fields $\BO_0/\bp_0=\BO/\bp=\CO/\fp=k$. In particular, $\CO_0/\fp_0=k$ and thus $\fp_0$ is ramified in $K/K_0$.

\noindent(3) By \cite[Chap. III, \S 6, Corollary 2]{Serre-localfields}, the different $\partial_{\CO/\BO}=(f'(\omega))$ is principal where $f(x)\in \BO[x]$ is the monic minimal polynomial of $\omega$.  Suppose 
\[t=\deg(f)=[K:\BK]=[\BQ(\omega):M]\geq 2.\]
 We compute that $f'(\omega)^{1-\dag}=-\N_{\BQ(\omega)/M}(\omega)\omega^{t-2}$. Let $\omega^{a_1}=\omega,\cdots, \omega^{a_t}$ be the Galois conjugates of $\omega$ over $\BK$. Since $\BK$ and $\BQ(\omega)$ are linearly disjoint over $M$, these are also Galois conjugates of $\omega$ over $M$. Then
\[f(x)=\prod_{i=1}^{t}(x-\omega^{a_i})=x^t+\cdots+(-1)^t\N_{\BQ(\omega)/M}(\omega),\]
\[f'(\omega)=\prod_{i=2}^{t}(\omega-\omega^{a_i}).\]
Suppose $g$ is the minimal polynomial of $\omega^{-1}$. Then
\[g(x)=\prod_{i=1}^{t}(x-\omega^{-a_i}),\quad g'(\omega^{-1})=\prod_{i=2}^{t}(\omega^{-1}-\omega^{-a_i})=f'(\omega)^\dag.\]
On the other hand
\[g(x)=(-1)^t\N_{\BQ(\omega)/M}(\omega)^{-1}x^tf(1/x).\]
Taking the derivative and substituting $x=\omega^{-1}$,
\[g'(\omega^{-1})=(-1)^{t-1}\N_{\BQ(\omega)/M}(\omega)^{-1}\omega^{-(t-2)}f'(\omega).\]
Thus
\[f'(\omega)/f'(\omega)^\dag=(-1)^{t-1}\N_{\BQ(\omega)/M}(\omega)\omega^{t-2}.\]
If $p=2$, then $t$ is even, and all of $-1$, $\N_{\BQ(\omega)/M}(\omega)$ and $\omega^{t-2}$ are squares and we denote 
\[(-1)^{t-1}\N_{\BQ(\omega)/M}(\omega)\omega^{t-2}=\omega^{2m}.\]
If $p\neq 2$, then $\N_{\BQ(\omega)/M}(\omega)\omega^{t-2}=\omega^{2m}$ is a square. In either case taking $\alpha=f'(\omega)\omega^{-m}$ will suffice our purpose.
\end{proof}

\subsection{Serre-tensor construction}\label{stc}
Assume the data $(A,\lambda,\BK,\BO,p,\bp)$  satisfies (\ref{CM}) with $\BO=\CO_\BK$ and $({\bf H1})$-$({\bf H3})$ but that $\omega\notin \BO$.  In certain cases we can use Serre-tensor construction to construct a new data $(A_1, \lambda_1, K, \CO, p,\fp)$ satisfying all these conditions (\ref{CM}) and $({\bf H1})$-$({\bf H3})$  with $\omega\in \CO$ as follows. 

As in the last subsection, let $K$ be the CM field $\BK(\omega)$ and $\CO=\BO[\omega]$. Then $\fp=(\bp,1-\omega)$ is a maximal ideal of $\CO$ with residue field $\CO/\fp=\BO/\bp=k$. We denote again the complex conjugation on $K$ by $\dag$, $K_0$ the maximal real subfield of $K$ and $M=\BK\cap \BQ(\omega)$. In the following we assume $p$ is prime to the discriminant $d_{\BK/M}$ and  if $p\neq 2$, $[\BQ(\omega):M]$ is odd.

Since $\BO=\CO_\BK$ is integrally closed, $\CO=\BO[\omega]$ is a finitely generated free $\BO$-module.  The functor $T\mapsto \CO \otimes_\BO A(T)$ for any $F$-schemes $T$ is representable by an abelian variety which we denote by $A_1=\CO\otimes_\BO A$. 
\begin{prop}
\begin{itemize}
\item[(1)] The abelian variety $A_1$ has CM by $K$ over $F$ and $\CO=K\cap \End_F(A_1)$. 
\item[(2)] We have $A_1[\fp]\cong A[\bp]$ as $G_F$-modules and $\Sel_{\fp}(A_1)\cong \Sel_\bp(A)$.
\end{itemize}
\end{prop}
\begin{proof}
By Theorem \ref{twist-av}, $A_1$ has CM by $K$ over $F$ with $\CO\hookrightarrow K\cap \End_F(A_1)$. By Lemma \ref{cyc-ext2}, $\CO$ is the maximal order of $K$, and thus $\CO=K\cap \End_F(A_1)$.  Also by Theorem \ref{twist-av}, $T_{\fp} A_1\cong \CO_\fp\otimes_{\BO_\bp} T_\bp(A)$, and thus
\[A_1[\fp]\cong \CO/\fp\otimes_\BO T_\bp(A)\cong \BO/\bp\otimes_\BO T_\bp(A)\cong A[\bp].\]
Through this isomorphism $A_1[\fp]\cong A[\bp]$, we have an identification $\RH^1(F,A_1[\fp])\cong \RH^1(F,A[\bp])$. From the commutative diagram
\[
\xymatrix{0\ar[r]&A[\bp]\ar[r]\ar[d]^{\rotatebox[origin=cc] {-90}{$\cong$}}&A\ar[d]\ar[r]&\bp^{-1}\otimes_\BO A\ar[d]\ar[r]&0\\
0\ar[r]&A_1[\fp]\ar[r]&A_1\ar[r]&{\fp}^{-1}\otimes_\CO A\ar[r]&0,}
\]
we see the local $\bp$-Kummer images coincide with the local $\fp$-Kummer images because $\CO/\fp\otimes_{\CO} A_1(F_v)\cong \BO/\bp\otimes_\CO A(F_v)$, and thus $\Sel_{\fp}(A_1)\cong \Sel_\bp(A)$.
\end{proof}

Fix $\partial_{\CO/\BO}=(\alpha)$ as in Lemma \ref{cyc-ext2} and let $h:  \CO \xrightarrow{\cong} \partial^{-1}_{\CO/\BO}=\Hom_\BO(\CO,\BO)$ be the non-degenerate hermitian form $h(x,y)=\Tr_{K/\BK}(xy^\dag/\alpha)$. Then $h$ is $\dag$-sesquilinear. Here the $\CO$-action of $\Hom_{\BO}(\CO,\BO)$ is given as $(af)(m)=f(am)$.  
\begin{prop}\label{symmetric-isogeny}
The tensor product $\lambda_1=h\otimes\lambda$ is a $\dag$-sesquilinear symmetric isogeny with kernel $\CO\otimes_\BO A[\lambda]$. In particular, $p\nmid \deg(\lambda_1)$.
\end{prop}
\begin{proof}
Since both $h$ and $\lambda$ are $\dag$-sesquilinear, so is $\lambda_1$. By assumption that $[\BQ(\omega):M]$ is odd if $p\neq 2$ and Lemma \ref{cyc-ext2}, $h$ is $+1$-hermitian.  By \cite[Proposition 7]{Amir18}, noting the normalizations of $\CO$-actions, $\lambda_1$ is a symmetric isogeny. Since $h:\CO \xrightarrow{\cong} \partial^{-1}_{\CO/\BO}$ is a $\dag$-sesquilinear isomorphism of $\CO$-modules and $\CO$ is free over $\BO$, $\Ker(\lambda_1)\cong\CO\otimes_\BO A[\lambda]$.
\end{proof}

Thus $(A_1, \lambda_1, K, \CO, p,\fp)$ satisfies (\ref{CM}) and $({\bf H1})$-$({\bf H3})$ as desired.

\bibliographystyle{plain}
\bibliography{reference}
\end{document}